%% file: IvCatArx.tex
\documentclass[11pt]{article}
\usepackage{amssymb,amsfonts,amsmath,graphicx,hyperref}

\input{ivdiag}
\input{tcilatex}
\begin{document}

\title{Applied Categories and Functors for Undergraduates}
\author{Vladimir G. Ivancevic\thanks{
Vladimir.Ivancevic@dsto.defence.gov.au} \and Tijana T. Ivancevic\thanks{ Tijana.Ivancevic@alumni.adelaide.edu.au}}
\date{}
\maketitle

\begin{abstract}
These are lecture notes for a 1--semester undergraduate course (in
computer science, mathematics, physics, engineering, chemistry or biology) in applied categorical meta-language. The only
necessary background for comprehensive reading of these notes are first-year calculus and linear algebra.
\end{abstract}

\tableofcontents

\bigbreak\bigbreak\bigbreak
\section{Introduction}
In modern mathematical sciences whenever one defines a new class
of mathematical objects, one proceeds almost in the next breath to
say what kinds of maps between objects will be considered
\cite{Eilenberg,MacLane,Switzer,GaneshSprBig,GaneshADG}.
A general framework for dealing with
situations where we have some \emph{objects} and \emph{maps
between objects}, like sets and functions, vector spaces and
linear operators, points in a space and paths between points, etc.
-- gives the modern metalanguage of categories and functors.
Categories are mathematical universes and functors are
`projectors' from one universe onto another.

\section{Sets and Maps}

\subsection{Notes from Set Theory}

Given a map (or, a function) $f:A\rightarrow B$, the set $A$ is
called the \emph{domain} of $f$, and denoted $\limfunc{Dom}f$. The
set $B$ is called the \emph{codomain} of $f$, and denoted
$\limfunc{Cod}f.$ The codomain is not to be confused with the
\emph{range} of $f(A)$, which is in general only a subset of $B$.

A map $f:X\rightarrow Y$ is called \emph{injective}, or 1--1, or an
\emph{injection}, iff for every $y$ in the codomain $Y$ there is \emph{at
most} one $x$ in the domain $X$ with $f(x)=y$. Put another way,
given $x$ and $x'$ in $X$, if $f(x)=f(x')$, then it follows that
$x=x'$. A map $f:X\rightarrow Y$ is called
\emph{surjective}, or \emph{onto}, or a \emph{surjection}, iff for every $%
y $ in the codomain $\limfunc{Cod}f$ there is \emph{at least} one $x$ in
the \emph{domain} $X$ with $f(x)=y$. Put another way, the
\emph{range} $f(X)$ is equal to the codomain $Y$. A map is
\emph{bijective} iff it is both injective and surjective.
Injective functions are called \emph{monomorphisms},
and surjective functions are called \emph{epimorphisms} in the
\emph{category of sets} (see below). Bijective functions are called \emph{isomorphisms}.

A \emph{relation} is any subset of a \emph{Cartesian product} (see
below). By definition, an \emph{equivalence relation} $\alpha$ on
a set $X$ is a relation which is \emph{reflexive, symmetrical} and
\emph{transitive}, i.e., relation that satisfies the following
three conditions:

\begin{enumerate}
  \item \emph{Reflexivity}: each element $x\in X$ is equivalent
  to itself, i.e.$,$ $x\alpha x$;
  \item \emph{Symmetry}: for any two elements $a,b\in X$, $a\alpha
b$ implies $b\alpha a$; ~and
  \item \emph{Transitivity}: $a\alpha b$\ and $b\alpha c$ implies $a\alpha c$.
\end{enumerate}

Similarly, a relation $\leq $ defines a \emph{partial order} on a
set $S$ if it has the following properties:

\begin{enumerate}
\item  \emph{Reflexivity}: $a\leq a$\ for all $a\in S$;

\item  \emph{Antisymmetry}: $a\leq b$\ and $b\leq a$ implies
$a=b$; ~and

\item  \emph{Transitivity}: $a\leq b$\ and $b\leq c$ implies
$a\leq c$.
\end{enumerate}

A \emph{partially ordered set} (or \emph{poset}) is a set taken
together with a partial order on it. Formally, a partially ordered
set is defined as an ordered pair $P=(X,\leq )$, where $X$ is
called the \emph{ground set} of $P$ and $\leq $\ is the partial
order of $P$.

\subsection{Notes From Calculus}

\subsubsection{Maps}

Recall that a \emph{map} (or, \emph{function}) $f$ is a
\emph{rule} that assigns to each element $x$ in a set $A$ exactly
one element, called $f(x)$, in a set $B$. A map could be thought
of as a \emph{machine} $[[f]]$ with $x-$input (the \emph{domain}
of $f$ is the set of all possible inputs) and $f(x)-$output (the
\emph{range} of $f$ is the set of all possible outputs)
\cite{Calc}
\begin{equation*}
x\rightarrow [[f]]\rightarrow f(x).
\end{equation*}
There are four possible ways to represent a function (or map): (i)
verbally (by a description in words); (ii) numerically (by a table
of values); (iii) visually (by a graph); and (iv) algebraically
(by an explicit formula). The most common method for visualizing a
function is its \emph{graph}. If $f$ is a function with domain
$A$, then its graph is the set of ordered input--output pairs
\begin{equation*}
\{(x,f(x)):x\in A\}.
\end{equation*}A generalization of the graph concept is a concept
of a \emph{cross--section of a fibre bundle}, which is one of the
core geometrical objects for dynamics of complex systems (see \cite{GaneshSprBig}).

\subsubsection{Algebra of Maps}

Let $f$ and $g$ be maps with domains $A$ and $B$. Then the maps
$f+g$, $f-g$, $fg$, and $f/g$ are defined as follows \cite{Calc}
\begin{eqnarray*}
(f+g)(x) &=&f(x)+g(x)\text{ \ \ \ \ \ \ \ \ \ \ \ domain }=A\cap B, \\
(f-g)(x) &=&f(x)-g(x)\text{ \ \ \ \ \ \ \ \ \ \ \ domain }=A\cap B, \\
(fg)(x) &=&f(x)\,g(x)\text{ \ \ \ \ \ \ \ \ \ \ \ domain }=A\cap B, \\
\left( \frac{f}{g}\right) (x) &=&\frac{f(x)}{g(x)}\text{ \ \ \ \ \
\ \ \ \ \ \ domain }=\{x\in A\cap B:g(x)\neq 0\}.
\end{eqnarray*}

\subsubsection{Compositions of Maps}

Given two maps $f$ and $g$, the composite map $f\circ g$,
called the \emph{composition} of $f$ and $g$, is defined by
\begin{equation*}
(f\circ g)(x)=f(g(x)).
\end{equation*}
The $(f\circ g)-$machine is composed of the $g-$machine (first) and then the $%
f-$machine \cite{Calc},
\begin{equation*}
x\rightarrow [[g]]\rightarrow g(x)\rightarrow [[f]]\rightarrow
f(g(x)).
\end{equation*}
For example, suppose that $y=f(u)=\sqrt{u}$ and $u=g(x)=x^{2}+1$.
Since $y$ is a function of $u$ and $u$ is a function of $x$, it
follows that $y$ is ultimately a function of $x$. We calculate
this by substitution
\begin{equation*}
y=f(u)=f\circ g=f(g(x))=f(x^{2}+1)=\sqrt{x^{2}+1}.
\end{equation*}

\subsubsection{The Chain Rule}

If $f$ and $g$ are both differentiable (or smooth, i.e.,
$C^{\infty}$) maps and $h=f\circ g$ is the composite map defined
by $h(x)=f(g(x))$, then $h$ is differentiable and $h'$ is given by
the product \cite{Calc}
\begin{equation*}
h'(x)=f'(g(x))\,g'(x).
\end{equation*}
In Leibniz notation, if $y=f(u)$ and $u=g(x)$ are both
differentiable maps, then
\begin{equation*}
\frac{dy}{dx}=\frac{dy}{du}\frac{du}{dx}.
\end{equation*}
The reason for the name \emph{chain rule} becomes clear if we add
another link to the chain. Suppose that we have one more
differentiable map $x=h(t)$. Then, to calculate the derivative of
$y$ with respect to $t$, we use the chain rule twice,
\begin{equation*}
\frac{dy}{dt}=\frac{dy}{du}\frac{du}{dx}\frac{dx}{dt}.
\end{equation*}

\subsubsection{Integration and Change of Variables}

Given a 1--1 continuous (i.e., $C^{0}$) map $F$ with a nonzero \emph{Jacobian} $\left| \frac{\partial (x,...)}{%
\partial (u,...)}\right|$\ that maps a region $S$ onto a region $R$ (see
\cite{Calc}), we have the following substitution formulas:\\
\smallskip

\noindent 1. For a single integral,
\begin{equation*}
\int_{R}f(x)\,dx=\int_{S}f(x(u))\frac{\partial x}{\partial u}du;
\end{equation*}
2. For a double integral,
\begin{equation*}
\iint_{R}f(x,y)\,dA=\iint_{S}f(x(u,v),y(u,v))\left| \frac{\partial (x,y)}{%
\partial (u,v)}\right| dudv;
\end{equation*}
3. For a triple integral,
\begin{equation*}
\iiint_{R}f(x,y,z)\,dV=\iiint_{S}f(x(u,v,w),y(u,v,w),z(u,v,w))\left| \frac{%
\partial (x,y,z)}{\partial (u,v,w)}\right| dudvdw;
\end{equation*}
4. Generalization to $n-$tuple integrals is obvious.

\subsection{Notes from General Topology}

\emph{Topology} is a kind of \emph{abstraction} of Euclidean
geometry, and also a natural framework for the study of
\emph{continuity}.\footnote{Intuitively speaking, a function $f:
\Bbb R \to \Bbb R $ is continuous near a point $x$ in its domain
if its value does not jump there. That is, if we just take $\delta
x$ to be small enough, the two function values $f(x)$ and $f(x +
\delta x)$ should approach each other arbitrarily closely. In more
rigorous terms, this leads to the following definition: A function
$f: \Bbb R  \to \Bbb R $ is continuous at $x \in \Bbb R $ if for
all $\epsilon > 0$, there exists a $\delta > 0$ such that for all
$y \in \Bbb R $ with $|y-x|<\delta$, we have that $|f(y) -
f(x)|<\epsilon$. The whole function is called continuous if it is
continuous at every point $x$.} Euclidean geometry is abstracted
by regarding triangles, circles, and squares as being the same
basic object. Continuity enters because in saying this one has in
mind a \emph{continuous deformation} of a triangle into a square
or a circle, or any arbitrary shape. On the other hand, a disk
with a hole in the center is topologically different from a circle
or a square because one cannot create or destroy holes by
continuous deformations. Thus using topological methods one does
not expect to be able to identify a geometrical figure as being a
triangle or a square. However, one does expect to be able to
detect the presence of gross features such as holes or the fact
that the figure is made up of two disjoint pieces etc. In this way
topology produces theorems that are usually qualitative in nature
-- they may assert, for example, the existence or non--existence
of an object. They will not, in general, give the means for its
construction \cite{Nash}.

\subsubsection{Topological Space}

Study of topology starts with the fundamental notion of
\emph{topological space}. Let $X$ be any \emph{set} and
$Y=\{X_{\alpha }\}$ denote a collection, finite or infinite of
subsets of $X$. Then $X$ and $Y$ form a topological space provided
the $X_{\alpha }$ and $Y$ satisfy:
\begin{enumerate}
\item  Any finite or infinite subcollection $\left\{ Z_{\alpha
}\right\} \subset X_{\alpha }$ has the property that $\cup
Z_{\alpha }\in Y$;

\item  Any \emph{finite subcollection} $\left\{ {Z_{\alpha
_{1}},...,Z_{\alpha _{n}}}\right\} \subset X_{\alpha }$ has the property that%
$\cap Z_{\alpha _{i}}\in Y$; ~and

\item  Both $X$ and the empty set belong to $Y$.
\end{enumerate}

The set $X$ is then called a topological space and the $X_{\alpha
}$ are called \emph{open sets}. The choice of $Y$ satisfying (2)
is said to give a topology to $X.$

Given two topological spaces $X$ and $Y$, a map $f:X\rightarrow Y$ is \emph{continuous} if the
inverse image of an open set in $Y$ is an open set in $X$.

The main general idea in topology is to study spaces which can be
continuously deformed into one another, namely the idea of
\emph{homeomorphism}. If we have two topological spaces $X$ and
$Y$, then a map $f:X\rightarrow Y$ is called a homeomorphism iff
\begin{enumerate}
    \item $f$ is continuous ($C^0$), ~and
    \item There exists an inverse of $f$, denoted $f^{-1}$, which is also continuous.
\end{enumerate}
Definition (2) implies that if $f$ is a homeomorphism then so is
$f^{-1}$. Homeomorphism is the main topological example of
\emph{reflexive}, \emph{symmetrical} and \emph{transitive
relation}, i.e., \emph{equivalence relation}. Homeomorphism
divides all topological spaces up into
\emph{equivalence classes}. In other words, a pair of topological spaces, $%
X$ and $Y$, belong to the same equivalence class if they are
homeomorphic.

The second example of topological equivalence relation is
\emph{homotopy}. While homeomorphism generates equivalence classes
whose members are topological spaces, homotopy generates
equivalence classes whose members are continuous ($C^0$) maps.
Consider two continuous maps $f,g:X\rightarrow Y$ between
topological spaces $X$ and $Y$. Then the map $f$ is said to be
\emph{homotopic} to the map $g$ if $f$ can be continuously
deformed into $g$ (see below for the precise definition of
homotopy). Homotopy is an equivalence relation which divides the
space of continuous maps between two topological spaces into
equivalence classes \cite{Nash}.

Another important notions in topology are \emph{covering},
\emph{compactness} and \emph{connectedness}. Given a family of
sets $\{X_{\alpha }\}=X$ say, then $X$ is a \emph{covering} of
another set $Y$ if $\cup X_{\alpha }$ contains $Y$. If all the
$X_{\alpha }$ happen to be open sets the covering is called an
\emph{open covering}. Now consider the set $Y$ and all its
possible open coverings. The set $Y$ is \emph{compact} if for
every open covering $\{X_{\alpha }\}$ with $\cup X_{\alpha
}\supset Y$ there
always exists a finite subcovering $\{X_{1},...,X_{n}\}$ of $Y$ with $%
X_{1}\cup ...\cup X_{n}\supset Y$. Again, we define a set $Z$ to
be \emph{connected} if it cannot be written as $Z=Z_{1}\cup
Z_{2}$, where $Z_{1}\ $and $Z_{2}$ are both open non--empty sets and $Z_{1}\cap
Z_{2}$ is an empty set.

Let $A_{1},A_{2},...,A_{n}$ be closed subspaces of a topological
space $X$\ such that $X=\cup _{i=1}^{n}A_{i}$. Suppose
$f_{i}:A_{i}\rightarrow Y$ is a function, $1\leq i\leq n$, such that
\begin{equation}
f_{i}|A_{i}\cap A_{j}=f_{j}|A_{i}\cap A_{j},\qquad (1\leq i,j\leq n).
\label{cont_over}
\end{equation}%
In this case $f$ is continuous iff each $f_{i}$ is. Using this
procedure we
can define a $C^0-$function $f:X\rightarrow Y$ by cutting up the space $%
X$ into closed subsets $A_{i}$ and defining $f$ on each $A_{i}$
separately in such a way that $f|A_{i}$ is obviously continuous;
we then have only to check that the different definitions agree on
the \emph{overlaps} $A_{i}\cap A_{j}$.

The \emph{universal property of the Cartesian product}: let $p
_{X}:X\times Y\rightarrow X$, and $p _{Y}:X\times Y\rightarrow Y$
be the \emph{projections} onto the first and second factors,
respectively. Given any pair of functions $f:Z\rightarrow X$ and
$g:Z\rightarrow Y$ there is a
unique function $h:Z\rightarrow X\times Y$ such that $p _{X}\circ h=f$, and $%
p _{Y}\circ h=g$. Function $h$ is continuous iff both $f$ and $g$
are. This property characterizes $X\times Y$ up to isomorphism.
In particular, to check that a given function $h:Z\rightarrow X$
is continuous it will suffice to check that $p _{X}\circ h$ and $p
_{Y}\circ h$ are continuous.

The \emph{universal property of the quotient}: let $\alpha $ be an
equivalence relation on a topological space $X$, let $X/\alpha $
denote the \emph{space of equivalence classes} and $p _{\alpha
}:X\rightarrow X/\alpha $ the \emph{natural projection}. Given a
function $f:X\rightarrow Y$, there is a function $f':X/\alpha
\rightarrow Y$ with $f'\circ p _{\alpha }=f$ iff $x\alpha x'$
implies $f(x)=f(x') $, for all $x\in X$. In this case $f'$ is
continuous iff $f$ is. This property characterizes $X/\alpha $ up
to homeomorphism.

\subsubsection{Homotopy}

Now we return to the fundamental notion of homotopy. Let $I$ be a
compact unit interval $I=[0,1]$. A \emph{homotopy} from $X$ to $Y$
is a continuous
function $F:X\times I\rightarrow Y$. For each $t\in I$ one has $%
F_{t}:X\rightarrow Y$ defined by $F_{t}(x)=F(x,t)$ for all $x\in
X$. The
functions $F_{t}$ are called the `stages' of the homotopy. If $%
f,g:X\rightarrow Y$ are two continuous maps, we say $f$ is
homotopic to $g$, and write $f\simeq g$, if there is a homotopy
$F:X\times I\rightarrow Y$ such that $F_{0}=f$ and $F_{1}=g$. In
other words, $f$ can be continuously deformed into $g$ through the
stages $F_{t}$. If $A\subset X$ is a subspace, then $F$ is a
homotopy relative to $A$ if $F(a,t)=F(a,0)$, for all $a\in A,t\in
I$.

The homotopy relation $\simeq $ is an equivalence relation. To
prove that we have $f\simeq f$ is obvious; take $F(x,t=f(x)$, for
all $x\in X,\,t\in I$. If $f\simeq g$ and $F$ is a homotopy from
$f$ to $g$, then $G:X\times I\rightarrow Y$ defined by
$G(x,t)=F(x,1-t)$, is a homotopy from $g$ to $f$, i.e., $g\simeq
f$. If $f\simeq g$ with homotopy $F$ and $g\simeq f$ with homotopy
$G$, then $f\simeq h$ with homotopy $H$ defined by
\begin{equation*}
H(x,t)=\left\{
\begin{array}{c}
F(x,t),\qquad \ \ \ \ \ \ \ 0\leq t\leq 1/2 \\
G(x,2t-1),\qquad 1/2\leq t\leq 1
\end{array}
\right..
\end{equation*}
To show that $H$ is continuous we use the relation
(\ref{cont_over}).

In this way, the set of all $C^0-$functions $f:X\rightarrow Y$
between two topological spaces $X$ and $Y$, called the
\emph{function space} and denoted by $Y^{X}$, is partitioned into
equivalence classes under the relation $\simeq $. The equivalence
classes are called \emph{homotopy classes}, the homotopy class of
$f$ is denoted by $[f]$, and the set of all homotopy classes is
denoted by $[X;Y]$.

If $\alpha $ is an equivalence relation on a topological space $X$ and $%
F:X\times I\rightarrow Y$ is a homotopy such that each stage
$F_{t}$ factors through $X/\alpha $, i.e., $x\alpha x'$ implies
$F_{t}(x)=F_{t}(x')$, then $F$ induces a homotopy $F':(X/\alpha
)\times I\rightarrow Y$ such that $F'\circ (p_{\alpha }\times
1)=F$.

Homotopy theory has a range of applications of its own, outside
topology and geometry, as for example in proving Cauchy theorem in
complex variable theory, or in solving nonlinear equations of
artificial neural networks.

A \emph{pointed set} $(S,s_{0})$ is a set $S$ together with a
distinguished point $s_{0}\in S$. Similarly, a \emph{pointed
topological space} $(X,x_{0})$ is a space $X $ together with a distinguished point $%
x_{0}\in X$. When we are concerned with pointed spaces $(X,x_{0}),(Y,y_{0})$, etc, we always require that all functions $f:X\rightarrow Y$
shell preserve base points, i.e., $f(x_{0})=y_{0}$, and that all homotopies $%
F:X\times I\rightarrow Y$ be relative to the base point, i.e., $%
F(x_{0},t)=y_{0}$, for all $t\in I$. We denote the homotopy
classes of base point--preserving functions by $[X,x_{0};Y,y_{0}]$
(where homotopies are
relative to $x_{0}).$ $[X,x_{0};Y,y_{0}]$ is a pointed set with base point $%
f_{0}$, the constant function: $f_{0}(x)=y_{0}$, for all $x\in X$.

A \emph{path} $\gamma (t)$ from $x_{0}$ to $x_{1}$ in a topological space $%
X$ is a continuous map $\gamma :I\rightarrow X$ with $\gamma (0)=x_{0}$ and $%
\gamma (1)=x_{1}$. Thus $X^{I}$ is the space of all paths in $X$
with the
compact--open topology. We introduce a relation $\sim $ on $X$ by saying $%
x_{0}\sim x_{1}$ iff there is a path $\gamma :I\rightarrow X$ from
$x_{0}$ to $x_{1}$. Clearly, $\sim $ is an equivalence relation;
the set of equivalence classes is denoted by $\pi _{0}(X)$.
The elements of $\pi _{0}(X) $ are called the \emph{path
components}, or $0-$\emph{components} of $X$. If $\pi _{0}(X)$
contains just one element, then $X$ is called \emph{path
connected}, or $0-$\emph{connected}. A \emph{closed path}, or
\emph{loop} in $X$ at the point $x_{0}$ is a path $\gamma (t)$ for
which $\gamma (0)=\gamma (1)=x_{0}.$ The \emph{inverse loop}
$\gamma
^{-1}(t)$ based at $x_{0}\in X$ is defined by $\gamma ^{-1}(t)=\gamma (1-t)$%
, for $0\leq t\leq 1.$ The \emph{homotopy of loops} is the
particular case of the above defined homotopy of continuous maps.

If $(X,x_{0})$ is a pointed space, then we may regard $\pi
_{0}(X)$ as a pointed set with the $0-$component\ of $x_{0}$ as a
base point. We use the notation $\pi _{0}(X,x_{0})$ to denote $p
_{0}(X,x_{0})$ thought of as a pointed set. If $f:X\rightarrow Y$
is a map then $f$ sends $0-$components of $X$ into $0-$components
of $Y$ and hence defines a function $\pi _{0}(f):\pi
_{0}(X)\rightarrow \pi _{0}(Y)$. Similarly, a base point--preserving map $%
f:(X,x_{0})\rightarrow (Y,y_{0})$ induces a map of pointed sets
$\pi _{0}(f):\pi _{0}(X,x_{0})\rightarrow \pi _{0}(Y,y_{0})$. In
this way defined $\pi _{0}$ represents a `functor' from the
`category' of topological (point) spaces to the underlying
category of (point) sets (see the next subsection).

The \emph{fundamental group} (introduced by Poincar\'e), denoted $\pi_{1}(X)$, of a pointed space $(X,x_{0})$ is the group (see Appendix) formed by the equivalence classes of the set of all \emph{loops}, i.e., closed homotopies with initial and final points at a given base point $x_0$. The identity element of this group is the set of all paths homotopic to the degenerate path consisting of the point $x_0$.\footnote{The group product $f*g$ of loop $f$ and loop $g$ is given by the path of $f$ followed by the path of $g$. The identity element is represented by the constant path, and the inverse $f^{-1}$ of $f$ is given by traversing $f$ in the opposite direction. The fundamental group $\pi_{1}(X)$ is independent of the choice of base point $x_0$ because any loop through $x_0$ is homotopic to a loop through any other point $x_1$.} The fundamental group $\pi_{1}(X)$ only depends on the homotopy type of the space $X$, that is, fundamental groups of homeomorphic spaces are isomorphic.

Combinations of topology and calculus give differential topology and
differential geometry.

\subsection{Commutative Diagrams}

The \emph{category theory} (see below) was born with an
observation that many properties of mathematical systems can be
unified and simplified by a presentation with \emph{commutative
diagrams of arrows} \cite{Eilenberg,MacLane}. Each arrow $f:X\rightarrow Y$
represents a function (i.e., a map, transformation, operator);
that is, a source (domain) set $X$, a target (codomain) set $Y$,
and a rule $x\mapsto f(x)$ which assigns to each element $x\in X$
an element $f(x)\in Y$. A typical diagram of sets and functions is
\begin{equation*}
\qtriangle[X`Y`Z;f`h`g\qquad\qquad \text{or} \qquad\qquad]
\qtriangle[X`f(X)`g(f(X));f`h`g]
\end{equation*}
This diagram is \emph{commutative} iff $h=g\circ f$, where $g\circ
f$ is the usual composite function $g\circ f:X\rightarrow Z$,
defined by $x\mapsto g(f(x))$.

Similar commutative diagrams apply in other mathematical, physical
and computing contexts; e.g., in the `category' of all topological
spaces, the letters $X,Y,$ and $Z$ represent topological spaces
while $f,g,$ and $h$ stand for continuous maps. Again, in the
category of all groups, $X,Y,$ and $Z$ stand for groups, $f,g,$
and $h$ for homomorphisms.

\begin{figure}[htb]
\centerline{\includegraphics[width=5cm]{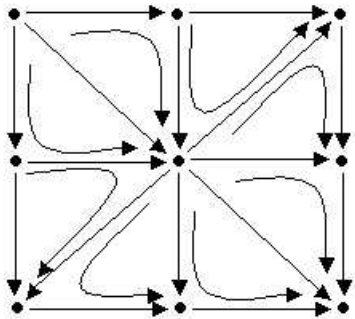}}
\caption{A
commutative flow (denoted by curved arrows) on a
triangulated digraph. Commutativity of the whole diagram
follows from commutativity of its triangular components.}
\label{ComDiag}
\end{figure}

Less formally, composing maps is like following directed paths
from one object to another (e.g., from set to set). In general, a
diagram is commutative iff any two paths along arrows that start
at the same point and finish at the same point yield the same
`homomorphism' via compositions along successive arrows.
Commutativity of the whole diagram follows from commutativity of
its triangular components (depicting a `commutative flow', see
Figure \ref{ComDiag}). Study of commutative diagrams is popularly
called `diagram chasing', and provides a powerful tool for
mathematical thought.

Many properties of mathematical constructions may be represented
by \emph{universal properties} of diagrams \cite{MacLane}.
Consider the \emph{Cartesian product} $X\times Y$ of two sets,
consisting as usual of all ordered pairs $\langle x,y\rangle$ of
elements $x\in X$ and $y\in Y$. The projections $\langle
x,y\rangle\mapsto x,\,\,\langle x,y\rangle\mapsto y$ of the
product on its `axes' $X$ and $Y$ are functions $p:X\times
Y\rightarrow X,\,\,q:X\times Y\rightarrow Y$. Any function
$h:W\rightarrow X\times Y$ from a third set $W$ is uniquely
determined by its composites $p\circ h$ and $q\circ h$.
Conversely, given $W$ and two functions $f$ and $g$ as in the
diagram below, there is a unique function $h$ which makes the
following diagram commute:
\begin{equation*}
\Atrianglepair<1`1`1`-1`1;>[W`X`X\times Y`Y;f`h`g`p`q]
\end{equation*}
This property describes the Cartesian product $X\times Y$
uniquely; the same diagram, read in the category of topological
spaces or of groups, describes uniquely the Cartesian product of
spaces or of the direct product of groups.

The construction `Cartesian product' is technically called a
`functor' because it applies suitably both to the sets and to the functions between them; two functions $%
k:X\rightarrow X'$ and $l:Y\rightarrow Y'$ have a function
$k\times l$ as their Cartesian product:
\begin{equation*}
k\times l:X\times Y\rightarrow X'\times Y',\qquad \langle
x,y\rangle \mapsto \langle kx,ly\rangle .
\end{equation*}

\section{Categories}
\label{category}

A category is a generic mathematical structure consisting of a
collection of \emph{objects} (sets with possibly additional
structure), with a corresponding collection of \emph{arrows}, or
\emph{morphisms}, between objects (agreeing with this additional
structure). A category $\mathcal{K}$
is defined as a pair $\left( \mathtt{Ob}(\mathcal{K}),\mathtt{Mor}(\mathcal{K})\right) $ of generic objects $A,B,\ldots $ in
$\mathtt{Ob}(\mathcal{K})$
and generic arrows $f:A\rightarrow B,\,g:B\rightarrow C,\ldots $ in $\mathtt{%
Mor}(\mathcal{K})$ between objects, with \textit{associative
composition}:
\begin{equation*}
A\cone{f}B\cone{g}C=A\cone{g\of f}C,
\end{equation*}
and \emph{identity} (\emph{loop}) arrow. (Note that in topological
literature, $\mathtt{Hom}(\mathcal{K})$ or
$\mathtt{hom}(\mathcal{K})$ is used instead of
$\mathtt{Mor}(\mathcal{K})$; see \cite{Switzer}).

A category $\mathcal{K}$ is usually depicted as a
\emph{commutative diagram} (i.e., a diagram with a common
\emph{initial object} $A$ and \emph{final object} $D$):\\
\smallskip
\begin{equation*}
\bfig \putsquare<1`1`1`1;500`400>(0,0)[A`B`C`D;f`h`g`k]
\put(250,200){\oval(800,700)}\put(230,170){$\mathcal{K}$}\efig
\end{equation*}

\bigbreak

\noindent To make this more precise, we say that a \emph{category}
$\mathcal{K}$ is defined if we have:

\begin{enumerate}
\item A \emph{class of objects} $\{A,B,C,...\}$ of $\mathcal{K}$,
denoted by $\mathtt{Ob}(\mathcal{K});$

\item A \emph{set of morphisms}, or \emph{arrows} $\mathtt{Mor}_{%
\mathcal{K}}(A,B),$ with elements $f:A\rightarrow B$, defined for
any \emph{ordered pair} $(A,B)\in \mathcal{K}$, such that for two
different
pairs $(A,B)\neq (C,D)$ in $\mathcal{K}$, we have $\mathtt{Mor}_{\mathcal{K}%
}(A,B)\cap \mathtt{Mor}_{\mathcal{K}}(C,D)=\emptyset $;

\item For any \emph{triplet} $(A,B,C)\in \mathcal{K}$ with
$f:A\rightarrow B$ and $g:B\rightarrow C$, there is a
\emph{composition} of morphisms
\begin{equation*}
\mathtt{Mor}_{\mathcal{K}}(B,C)\times
\mathtt{Mor}_{\mathcal{K}}(A,B)\ni (g,f)\rightarrow g\circ f\in
\mathtt{Mor}_{\mathcal{K}}(A,C),
\end{equation*}
written schematically as
\begin{equation*}
\frac{f:A\rightarrow B,\qquad g:B\rightarrow C}{g\circ
f:A\rightarrow C}.
\end{equation*}
\end{enumerate}

Recall from above that if we have a morphism $f\in \mathtt{Mor}_{\mathcal{K}}(A,B)$,
(otherwise written $f:A\rightarrow B$, or $A\cone{f}B$), then
$A=\limfunc{dom}(f)$ is a \emph{domain} of $f$, and
$B=\limfunc{cod}(f)$ is a \emph{codomain} of $f$ (of which
\emph{range} of $f$ is a subset, $B=\limfunc{ran}(f)$).
\smallskip

To make $\mathcal{K}$ a category, it must also fulfill the
following two properties:

\begin{enumerate}
\item \emph{Associativity of morphisms}\index{associativity of morphisms}: for all $f\in \mathtt{Mor}_{%
\mathcal{K}}(A,B)$, $g\in \mathtt{Mor}_{\mathcal{K}}(B,C)$, and
$h\in \mathtt{Mor}_{\mathcal{K}}(C,D)$, we have $h\circ (g\circ
f)=(h\circ g)\circ f$; in other words, the following diagram is
commutative
\begin{equation*}
\square <1`1`-1`1;1200`500>[A`D`B`C;h\circ (g\circ f)=(h\circ
g)\circ f`f`h`g]
\end{equation*}

\item \emph{Existence of identity morphism}: for every object
$A\in \mathtt{Ob}(\mathcal{K})$ exists a unique identity morphism
$1_{A}\in \mathtt{Mor}_{\mathcal{K}}(A,A)$; for any two
morphisms
$f\in \mathtt{Mor}_{\mathcal{K}}(A,B)$, and\\ $g\in \mathtt{Mor}_{\mathcal{K}%
}(B,C)$, compositions with identity morphism $1_{B}\in \mathtt{Mor}_{%
\mathcal{K}}(B,B)$ give $1_{B}\circ f=f$ and $g\circ 1_{B}=g$,
i.e., the following diagram is commutative:
\begin{equation*}
\Vtrianglepair<1`1`1`1`-1;>[A`B`C`B;f`g`f`1_{B}`g]
\end{equation*}
\end{enumerate}

The set of all morphisms of the category $\mathcal{K}$ is denoted
\begin{equation*}
\mathtt{Mor}(\mathcal{K})=\bigcup_{{A,B\in Ob(\mathcal{K})}}\mathtt{Mor}_{%
\mathcal{K}}(A,B).
\end{equation*}

If for two morphisms $f\in \mathtt{Mor}_{\mathcal{K}}(A,B)$ and
$g\in \mathtt{Mor}_{\mathcal{K}}(B,A)$ the equality $g\circ
f=1_{A}$ is valid, then the morphism $g$ is said to be \emph{left
inverse} (or \emph{retraction}), of $f$, and $f$ \emph{right
inverse} (or \emph{section}) of $g$. A morphism which is both
right and left inverse of $f$ is said to be \emph{two--sided
inverse} of $f$.

A morphism $m:A\rightarrow B$ is called \emph{monomorphism} in
$\mathcal{K}$ (i.e., \emph{1--1}, or \emph{injection} map), if for
any two parallel morphisms $f_1,f_2:C\rightarrow A$ in
$\mathcal{K}$ the equality $m\circ f_1=m\circ f_2$ implies
$f_1=f_2$; in other words, $m$ is monomorphism if it is \emph{left
cancellable}. Any morphism with a left inverse is monomorphism.

A morphism $e:A\rightarrow B$ is called \emph{epimorphism} in
$\mathcal{K}$ (i.e., \emph{onto}, or \emph{surjection} map), if
for any two morphisms $g_1,g_2:B\rightarrow C$ in $\mathcal{K}$
the equality $g_1\circ e=g_2\circ e$ implies $g_1=g_2$; in other
words, $e$ is epimorphism if it is \emph{right cancellable}. Any
morphism with a right inverse is epimorphism.

A morphism $f:A\rightarrow B$ is called \emph{isomorphism} in
$\mathcal{K}$ (denoted as $f:A\cong B$) if there exists a morphism
$f^{-1}:B\rightarrow A$ which is a two--sided inverse of $f$ in
$\mathcal{K}$. The relation of isomorphism is reflexive,
symmetric, and transitive, that is, an equivalence relation.

For example, an isomorphism in the category of sets is called a
set--isomorphism, or a \emph{bijection}, in the category of
topological spaces is called a topological isomorphism, or a
\emph{homeomorphism}, in the category of differentiable manifolds
is called a differentiable isomorphism, or a
\emph{diffeomorphism}.

A morphism $f\in \mathtt{Mor}_{\mathcal{K}}(A,B)$ is
\emph{regular} if there exists a morphism $g:B\rightarrow A$ in
$\mathcal{K}$ such that $f\circ g\circ f=f$. Any morphism with
either a left or a right inverse is regular.

An object $T$ is a \emph{terminal object} in $\mathcal{K}$ if to
each object $A\in \mathtt{Ob}(\mathcal{K})$ there is exactly one
arrow $A\rightarrow T$. An object $S$
is an \emph{initial object} in $\mathcal{K}$ if to each object $A\in \mathtt{Ob}(%
\mathcal{K})$ there is exactly one arrow $S\rightarrow A$. A \emph{null object} $%
Z\in \mathtt{Ob}(\mathcal{K})$ is an object which is both initial
and terminal; it is unique up to isomorphism. For any two objects
$A,B\in \mathtt{Ob}(\mathcal{K})$ there is a unique morphism
$A\rightarrow Z\rightarrow B$ (the composite through $Z$), called
the \emph{zero morphism} from $A$ to $B$.

A notion of subcategory is analogous to the notion of subset. A
subcategory $\mathcal{L}$ of a category $\mathcal{K}$ is said to
be a \emph{complete subcategory} iff for any objects $A,B\in
\mathcal{L}$, every morphism $A\rightarrow B$ of $\mathcal{L}$ is
in $\mathcal{K}$.

A \emph{groupoid} is a category in which every morphism is
invertible. A typical groupoid is the \emph{fundamental groupoid}
$\Pi_1(X)$ of a topological space $X$. An object of $\Pi_1(X)$ is
a point $x\in X$, and a morphism $x\rightarrow x'$ of $\Pi_1(X)$
is a homotopy class of paths $f$ from $x$ to $x'$. The
\emph{composition} of paths $g:x'\rightarrow x''$ and
$f:x\rightarrow x'$ is the path $h$ which is `$f$ followed by
$g$'. Composition applies also to homotopy classes, and makes
$\Pi_1(X)$ a category and a groupoid (the inverse of any path is
the same path traced in the opposite direction).

A \emph{group} is a groupoid with one object, i.e., a
\emph{category with one object} in which \emph{all morphisms are
isomorphisms} (see Appendix). Therefore, if we try to generalize the concept of a
group, keeping associativity as an essential property, we get the
notion of a category.

A category is \emph{discrete} if every morphism is an identity. A
\emph{monoid} is a category with one object, which is a group without inverses. A group is a category with one object in which every morphism has a two--sided
inverse under composition.

\emph{Homological algebra}\index{homological algebra} was the
progenitor of category theory (see e.g., \cite{Dieudonne2}).
Generalizing L. Euler's formula: $f+v=e+2$, for the faces $f$, vertices $v$
and edges $e$ of a convex polyhedron, E. Betti defined \emph{numerical
invariants of spaces} by formal addition and subtraction of faces
of various dimensions. H. Poincar\'{e} formalized these and
introduced the concept of \emph{homology}. E. Noether stressed the fact that these
calculations go on in Abelian groups, and that the operation
$\partial _{n}$ taking a face of dimension $n$ to the alternating
sum of faces of dimension $n-1$ which form its boundary is a
homomorphism, and it also satisfies \emph{the boundary of a boundary is zero} rule: $\partial _{n}\circ
\partial _{n+1}=0$. There are many ways of approximating a given space by
polyhedra, but the quotient $H_{n}=\limfunc{Ker}\partial _{n}/\limfunc{Im}%
\partial _{n+1}$ is an invariant, the \textit{homology group}.

As a physical example from \cite{Coecke1,Coecke2}, consider some physical system of type
$A$ (e.g., an electron) and perform some \emph{physical operation} $f$ on it (e.g., perform a measurement on it), which results in a possibly different system $B$ (e.g., a perturbed electron), thus having a map $f:A\to B$. In a same way, we can perform a consecutive operation $g:B\to C$ (e.g., perform the second measurement, this time on $B$), possibly resulting in a different system $C$ (e.g., a secondly perturbed electron). Thus, we have a composition: $k=g\circ f$, representing the consecutive application of these two physical operations, or the following diagram commutes:
\begin{equation*}
\qtriangle[A`B`C;f`k`g]
\end{equation*}
In a similar way, we can perform another consecutive operation $h:C\to D$ (e.g., perform the third measurement, this time on $C$), possibly resulting in a different system $D$ (e.g., a thirdly perturbed electron). Clearly we have an associative composition $(h\circ g)\circ f = h\circ (g\circ f)$, or the following diagram commutes:
\begin{equation*}
\square <1`1`-1`1;1200`500>[A`D`B`C;h\circ (g\circ f)=(h\circ
g)\circ f`f`h`g]
\end{equation*}
Finally, if we introduce a trivial operation $1_A\in \mathtt{Mor}_{\mathcal{K}}(A,A)$, meaning `doing nothing on a system of type $A$', we have $1_B\circ f = f\circ 1_A = f$. In this way, we have constructed a generic physical category (for more details, see \cite{Coecke1,Coecke2}).

For the same \emph{operational} reasons, categories could be expected to play an important role in other fields where operations/processes play a central role: e.g., Computer Science (computer programs as morphisms) and Logic \& Proof Theory (proofs as morphisms). In the theoretical
counterparts to these fields category theory has become quite common practice (see \cite{Abramsky}).

\section{Functors}
\label{functor}

In algebraic topology, one attempts to assign to every topological
space $X$ some algebraic object $\mathcal{F}(X)$ in such a way
that to every
$C^0-$function $f:X\rightarrow Y$ there is assigned a homomorphism $\mathcal{F}%
(f):\mathcal{F}(X)\to\mathcal{F}(Y)$ (see
\cite{Switzer,GaneshSprBig}). One advantage of this procedure is,
e.g., that if one is trying to prove the non--existence of a
$C^0-$function $f:X\rightarrow Y$ with certain properties, one may
find it relatively easy to prove the non--existence of the
corresponding algebraic function $\mathcal{F}(f)$ and hence deduce
that $f$ could not exist. In other words, $\mathcal{F}$ is to be a
`homomorphism'
from one category (e.g., $\mathcal{T}$) to another (e.g., $\mathcal{G}$ or $%
\mathcal{A}$). Formalization of this notion is a \emph{functor}.

A functor is a generic \emph{picture} projecting (all objects and
morphisms of) a source category into
a target category. Let $\mathcal{K}=\left( \mathtt{Ob}(\mathcal{K}),\mathtt{Mor}(%
\mathcal{K})\right) $ be a \emph{source} (or domain)
\emph{category} and $\mathcal{L}=\left(
\mathtt{Ob}(\mathcal{L}),\mathtt{Mor}(\mathcal{L})\right)
$ be a \emph{target} (or codomain) category. A functor $\mathcal{F}=(%
\mathcal{F}_{O},\mathcal{F}_{M})$ is defined as a pair of maps, $\mathcal{F}%
_{O}:\mathtt{Ob}(\mathcal{K})\rightarrow \mathtt{Ob}(\mathcal{L})$ and $%
\mathcal{F}_{M}:\mathtt{Mor}(\mathcal{K})\rightarrow \mathtt{Mor}(\mathcal{L})$%
, preserving categorical symmetry (i.e., commutativity of all diagrams) of $%
\mathcal{K}$ in $\mathcal{L}$.

More precisely, a \emph{covariant functor}, or simply a \emph{functor}, $%
\mathcal{F}_{\ast }:\mathcal{K}\rightarrow \mathcal{L}$ is a
\emph{picture} in the target category $\mathcal{L}$ of (all
objects and morphisms of) the source category $\mathcal{K}$:
\smallskip
\begin{equation*}
\bfig \putsquare<1`1`1`1;500`400>(0,0)[A`B`C`D;f`h`g`k]
\put(250,200){\oval(800,700)}\put(230,170){$\mathcal{K}$} \putsquare%
<1`1`1`1;800`400>(1600,0)[\mathcal{F}(A)`\mathcal{F}(B)`\mathcal{F}(C)`
\mathcal{F}(D);\mathcal{F}(f)`\mathcal{F}(h)`\mathcal{F}(g)`\mathcal{F}(k)]
\put(2000,200){\oval(1350,700)}\put(1980,180){$\mathcal{L}$} \put(800,170){$%
\cone{\mathcal{F}_\ast}$} \efig
\end{equation*}
\bigbreak\smallskip

Similarly, a \emph{contravariant functor}, or a \emph{cofunctor}, $%
\mathcal{F}^{\ast }:\mathcal{K}\rightarrow \mathcal{L}$ is a
\emph{dual picture} with reversed arrows:\\
\begin{equation*}
\bfig \putsquare<1`1`1`1;500`400>(0,0)[A`B`C`D;f`h`g`k]
\put(250,200){\oval(800,700)}\put(230,170){$\mathcal{K}$} \putsquare%
<-1`-1`-1`-1;800`400>(1600,0)[\mathcal{F}(A)`\mathcal{F}(B)`\mathcal{F}(C)`
\mathcal{F}(D);\mathcal{F}(f)`\mathcal{F}(h)`\mathcal{F}(g)`\mathcal{F}(k)]
\put(2000,200){\oval(1350,700)}\put(1980,180){$\mathcal{L}$} \put(800,170){$%
\cone{\mathcal{F}^\ast}$} \efig
\end{equation*}
\bigbreak

In other words, a \emph{functor}
$\mathcal{F}:\mathcal{K}\rightarrow \mathcal{L}$ from a
\emph{source} category $\mathcal{K}$ to a \emph{target} category
$\mathcal{L}$, is a pair $\mathcal{F}=(\mathcal{F}_{O},
\mathcal{F}_{M})$ of maps $\mathcal{F}_{O}:\mathtt{Ob}(\mathcal{K}%
)\rightarrow \mathtt{Ob}(\mathcal{L})$, $\mathcal{F}_{M}:\mathtt{Mor}(%
\mathcal{K})\rightarrow \mathtt{Mor}(\mathcal{L})$, such that

\begin{enumerate}
\item If $f\in \mathtt{Mor}_{\mathcal{K}}(A,B)$ then
$\mathcal{F}_{M}(f)\in
\mathtt{Mor}_{\mathcal{L}}(\mathcal{F}_{O}(A),\mathcal{F}_{O}(B))$
in case
of the \emph{covariant} functor $\mathcal{F}_{\ast }$, and $\mathcal{F}%
_{M}(f)\in
\mathtt{Mor}_{\mathcal{L}}(\mathcal{F}_{O}(B),\mathcal{F}_{O}(A))$
in case of the \emph{contravariant} functor $\mathcal{F}^{\ast }$;

\item For all $A\in \mathtt{Ob}(\mathcal{K}):$ $\mathcal{F}_{M}(1_{A})=1_{%
\mathcal{F}_{O}(A)};$

\item For all $f,g\in \mathtt{Mor}(\mathcal{K})$: if $\limfunc{cod}(f)=\limfunc{dom%
}(g)$, then $\mathcal{F}_{M}(g\circ
f)=\mathcal{F}_{M}(g)\circ \mathcal{F}_{M}(f)$ in
case of the \emph{covariant} functor $\mathcal{F}_{\ast }$, and $\mathcal{F%
}_{M}(g\circ f)=\mathcal{F}_{M}(f)\circ \mathcal{F}_{M}(g)$ in
case of the \emph{contravariant} functor $\mathcal{F}^{\ast }$.
\end{enumerate}

Category theory originated in algebraic topology, which tried to
assign algebraic invariants to topological structures. The golden
rule of such \emph{invariants} is that they should be
\emph{functors}. For example, the \emph{fundamental group} $\pi
_{1}$ is a functor. Algebraic topology constructs a group called
the \emph{fundamental group} $\pi _{1}(X)$ from any topological
space $X$, which keeps track of how many holes the space $X$
has. But also, any map between topological spaces determines a homomorphism $%
\phi :\pi _{1}(X)\rightarrow \pi _{1}(Y)$ of the fundamental
groups. So the fundamental group is really a functor $\pi
_{1}:\mathcal{T}\rightarrow \mathcal{G}$. This allows us to
completely transpose any situation involving \emph{spaces} and
\emph{continuous maps} between them to a parallel situation
involving \emph{groups} and \emph{homomorphisms} between them, and
thus reduce some topology problems to algebra problems.

Also, singular homology in a given dimension $n$ assigns to each
topological space $X$ an Abelian group $H_{n}(X)$, its $n$th
\emph{homology group} of $X$, and also to each continuous map
$f:X\rightarrow Y$ of spaces a corresponding homomorphism
$H_{n}(f):H_{n}(X)\rightarrow H_{n}(Y)$ of
groups, and this in such a way that $H_{n}(X)$ becomes a functor $H_{n}:%
\mathcal{T}\rightarrow \mathcal{A}$.

The leading idea in the \emph{use of functors in topology} is that
$H_n$ or $\pi_n$ gives an algebraic picture or image not just of
the topological spaces $X,Y$ but also of all the continuous maps
$f:X\rightarrow Y$ between them.

Similarly, there is a functor $\Pi _{1}:\mathcal{T}\rightarrow
\mathcal{G}$, called the `fundamental groupoid functor', which
plays a very basic role in algebraic topology. Here's how we get
from any space $X$ its `fundamental groupoid' $\Pi _{1}(X)$. To
say what the groupoid $\Pi _{1}(X)$ is, we need to say what its
objects and morphisms are. The objects in $\Pi _{1}(X)$ are just
the \emph{points} of $X$ and the morphisms are just certain
equivalence classes of \emph{paths} in $X$. More precisely, a morphism $%
f:x\rightarrow y$ in $\Pi _{1}(X)$ is just an equivalence class of
continuous paths from $x$ to $y$, where two paths from $x$ to $y$
are decreed equivalent if one can be continuously deformed to the
other while not moving the endpoints. (If this equivalence
relation holds, we say the two paths are `homotopic', and we call
the equivalence classes `homotopy classes of paths'; see
\cite{MacLane,Switzer}).

Another examples are covariant \emph{forgetful} functors:

\begin{itemize}
\item From the category of topological spaces to the category of sets;
it `forgets' the topology--structure.

\item From the category of metric spaces to the category of
topological spaces with the topology induced by the metrics; it
`forgets' the metric.
\end{itemize}

For each category $\mathcal{K}$, the \emph{identity functor} $I_{\mathcal{K%
}}$ takes every $\mathcal{K}-$object and every
$\mathcal{K}-$morphism to itself.

Given a category $\mathcal{K}$ and its subcategory $\mathcal{L}$,
we have an \emph{inclusion functor}
$\mathtt{In}:\mathcal{L}\to\mathcal{K}$.

Given a category $\mathcal{K}$, a \emph{diagonal functor} $\Delta:\mathcal{K}\to\mathcal{K}\times\mathcal{K}$ takes each object $A\in\mathcal{K}$ to the
object $(A,A)$ in the product category
$\mathcal{K}\times\mathcal{K}$.

Given a category $\mathcal{K}$ and a category of sets
$\mathcal{S}$, each
object $A\in\mathcal{K}$ determines a \emph{covariant Hom--functor} $%
\mathcal{K}[A,\_]:\mathcal{K}\rightarrow \mathcal{S}$, a
\emph{contravariant Hom--functor}
$\mathcal{K}[\_,A]:\mathcal{K}\to\mathcal{S}$, and a
\emph{Hom--bifunctor}
$\mathcal{K}[\_,\_]:\mathcal{K}^{op}\times\mathcal{K}\rightarrow
\mathcal{S}$.

A functor $\mathcal{F}:\mathcal{K}\rightarrow \mathcal{L}$ is a
\emph{faithful functor} if for all $A,B\in
\mathtt{Ob}(\mathcal{K})$ and for all $f,g\in
\mathtt{Mor}_{\mathcal{K}}(A,B)$, $\mathcal{F}(f)=\mathcal{F}(g)$ implies $%
f=g$; it is a \emph{full functor} if for every $h\in \mathtt{Mor}_{\mathcal{L}}(%
\mathcal{F}(A),\mathcal{F}(B))$, there is $g\in \mathtt{Mor}_{\mathcal{K}%
}(A,B)$ such that $h=\mathcal{F}(g)$; it is a \emph{full
embedding} if it is both full and faithful.

A \emph{representation of a group} is a functor
$\mathcal{F}:\mathcal{G} \rightarrow \mathcal{V}$. Thus, a category
is a generalization of a group and group representations are a
special case of category representations.

\section{Natural Transformations} \label{natural}

A \emph{natural transformation} (i.e., a \emph{functor morphism}) $\mathbf{\tau }%
: \mathcal{F}\overset{\cdot}{\rightarrow}\mathcal{G}$ is a
\emph{map
between two functors of the same variance}, $(\mathcal{F},\mathcal{G}):%
\mathcal{K}\rightrightarrows \mathcal{L}$, preserving categorical symmetry:%
\smallskip
\begin{equation*}
\bfig \putmorphism(0,200)(1,0)[A`B`f]{500}1a
\put(250,200){\oval(800,700)}\put(230,0){$\mathcal{K}$} \put(800,350){$\cone{%
\mathcal{F}}$} \put(930,200){$\mathbf{\tau}\Downarrow$} \put(800,-50){$\cone{%
\mathcal{G}}$} \putsquare<1`1`1`1;800`400>(1600,0)[\mathcal{F}(A)`\mathcal{F}%
(B)`\mathcal{G}(A)`\mathcal{G}(B);\mathcal{F}(f)`\mathbf{\tau }_{A}`\mathbf{%
\tau }_{B}`\mathcal{G}(f)] \put(2000,200){\oval(1200,700)}\put(1980,180){$%
\mathcal{L}$}\efig
\end{equation*}
\bigbreak

\noindent More precisely, all functors of the same variance from a source category $%
\mathcal{K}$ to a target category $\mathcal{L}$ form themselves
objects of
the \emph{functor category} $\mathcal{L}^{\mathcal{K}}$. Morphisms of $%
\mathcal{L}^{\mathcal{K}}$, called \emph{natural transformations},
are defined as follows.

Let $\mathcal{F}:\mathcal{K}\rightarrow \mathcal{L}$ and $\mathcal{G}:%
\mathcal{K}\rightarrow \mathcal{L}$ be two functors of the same
variance from a category $\mathcal{K}$ to a category
$\mathcal{L}$. Natural transformation
$\mathcal{F\overset{\mathbf{\tau }}{\longrightarrow }G}$ is a
family of morphisms such that for all $f\in
\mathtt{Mor}_{\mathcal{K}}(A,B) $
in the source category $\mathcal{K}$, we have $\mathcal{G}(f)\circ \mathbf{%
\tau }_{A}=\mathbf{\tau }_{B}\circ \mathcal{F}(f)$ in the target category $%
\mathcal{L}$. Then we say that the \emph{component} $\mathbf{\tau }_{A}:%
\mathcal{F}(A)\rightarrow \mathcal{G}(A)$ \emph{is natural in}
$A$.

If we think of a functor $\mathcal{F}$ as giving a \emph{picture}
in the target category $\mathcal{L}$ of (all the objects and
morphisms of) the source category $\mathcal{K}$, then a natural
transformation $\mathbf{\tau }$ represents a set of morphisms
mapping the picture $\mathcal{F}$ to another picture
$\mathcal{G}$, preserving the commutativity of all diagrams.

An invertible natural transformation, such that all components
$\mathbf{\tau }_{A}$ are isomorphisms) is called a \emph{natural
equivalence} (or, \emph{natural isomorphism}). In this case, the
inverses $(\mathbf{\tau }_{A})^{-1}$
in $\mathcal{L}$ are the components of a natural isomorphism $(\mathbf{\tau }%
)^{-1}:\mathcal{G}\overset{\ast }{\longrightarrow }\mathcal{F}$.
Natural equivalences are among the most important
\emph{metamathematical constructions} in algebraic topology (see
\cite{Switzer}).

As a mathematical example, let $\mathcal{B}$ be the category of Banach spaces over $%
\mathbb{R}$ and bounded linear maps. Define
$D:\mathcal{B}\rightarrow \mathcal{B} $ by taking $D(X)=X^{\ast
}=$ Banach space of bounded linear functionals on a space $X$ and
$D(f)=f^{\ast }$ for $f:X\rightarrow Y$ a bounded linear map. Then
$D$ is a cofunctor. $D^{2}=D\circ D$ is also a functor. We also
have the identity functor $1:\mathcal{B}\rightarrow \mathcal{B}$.
Define $T:1\rightarrow D\circ D$ as follows: for every $X\in
\mathcal{B}$ let $T(X):X\rightarrow D^{2}X=X^{\ast \ast }$ be the
\emph{natural inclusion} -- that is, for $x\in X$ we have
$[T(X)(x)](f)=f(x)$ for every $f\in X^{\ast }$. $T$ is a natural
transformation. On the subcategory of $n$D Banach spaces $T$ is
even a natural equivalence. The largest subcategory of
$\mathcal{B}$ on which $T$ is a natural equivalence is called the
category of reflexive Banach spaces \cite{Switzer}.

As a physical example, when we want to be able to conceive two physical
systems $A$ and $B$ as one whole (see \cite{Coecke1,Coecke2}), we can denote this using a (symmetric)
monoidal tensor product $A\otimes B$, and hence also
need to consider the compound operations
$$A\otimes B\cone{f\otimes g}C\otimes D,
$$
inherited from the operations on the individual systems.
Now, a (symmetric) \emph{monoidal category} is a category $\mathcal{K}$ defined as a pair $\left( \mathtt{Ob}(\mathcal{K}),\mathtt{Mor}(\mathcal{K})\right) $ of generic objects $A,B,\ldots $ in
$\mathtt{Ob}(\mathcal{K})$
and generic arrows $f:A\rightarrow B,\,g:B\rightarrow C,\ldots $ in $\mathtt{Mor}(\mathcal{K})$ between objects, defined using the symmetric monoidal tensor product:
\begin{eqnarray*}
\mathtt{Ob}(\mathcal{K})&:&\{A,B\} ~\mapsto~ A\otimes B,\\
\mathtt{Mor}(\mathcal{K})&:&\{A\cone{f}B, C\cone{g}D\} ~\mapsto~ A\otimes C\cone{f\otimes g}B\otimes D,
\end{eqnarray*}
with the additional notion of \emph{bifunctoriality}: if we apply an operation
$f$ to one system and an operation $g$ to another system, then the order in
which we apply them does not matter; that is, the following diagram commutes:
\begin{equation*}
\square <1`1`-1`1;1200`500>[A_1\otimes A_2`B_1\otimes A_2`A_1\otimes B_2`B_1\otimes B_2;f\otimes 1_{A_2}`1_{A_1}\otimes g`1_{B_1}\otimes g `f\otimes 1_{B_2}]
\end{equation*} which shows that both paths yield the same result (see \cite{Coecke1,Coecke2} for technical details).

As `categorical fathers', S. Eilenberg and S. MacLane, first observed, `category' has been
defined in order to define `functor' and `functor' has been
defined in order to define `natural transformations'
\cite{Eilenberg,MacLane}).

\subsection{Compositions of Natural Transformations}

Natural transformations can be \emph{composed} in two different
ways. First, we have an `ordinary' composition: if $\mathcal{F},\mathcal{G}$ and $%
\mathcal{H}$ are three functors from the source category
$\mathcal{A}\ $to
the target category $\mathcal{B}$, and then $\mathbf{\alpha }:\mathcal{F}%
\overset{\cdot }{\rightarrow }\mathcal{G}$, $\mathbf{\beta }:\mathcal{G}%
\overset{\cdot }{\rightarrow }\mathcal{H}$ are two natural
transformations, then the formula
\begin{equation}
\left( \beta \circ \alpha \right) _{A}=\beta _{A}\circ \alpha
_{A},\text{ \ \ \ \ \ \ (for all }A\in \mathcal{A}),
\label{natCom}
\end{equation}
defines a new natural transformation $\beta \circ \alpha :\mathcal{F}%
\overset{\cdot }{\rightarrow }\mathcal{H}$. This composition law
is clearly
associative and possesses a unit $1_{\mathcal{F}}$\ at each functor $%
\mathcal{F}$, whose $\mathcal{A}$--component is
$1_{\mathcal{FA}}.$

Second, we have the \emph{Godement product} of natural
transformations, usually denoted by $\ast $. Let $\mathcal{A}$, $\mathcal{B}$ and $%
\mathcal{C}$ be three categories, $\mathcal{F},\mathcal{G}$,
$\mathcal{H}$
and $\mathcal{K}$ be four functors such that $(\mathcal{F},\mathcal{G}):%
\mathcal{A}\rightrightarrows \mathcal{B}$ and $(\mathcal{H},\mathcal{K}):%
\mathcal{B}\rightrightarrows \mathcal{C}$, and $\mathbf{\alpha }:\mathcal{F}%
\overset{\cdot }{\rightarrow }\mathcal{G}$, $\mathbf{\beta }:\mathcal{H}%
\overset{\cdot }{\rightarrow }\mathcal{K}$ be two natural
transformations. Now, instead of (\ref{natCom}), the Godement
composition is given by
\begin{equation}
\left( \beta \ast \alpha \right) _{A}=\beta _{GA}\circ H\left(
\alpha _{A}\right) =K\left( \alpha _{A}\right) \circ \beta
_{FA},\text{ \ \ \ \ \ \ (for all }A\in \mathcal{A}),
\label{GodCom}
\end{equation}
which defines a new natural transformation $\beta \ast \alpha :\mathcal{H}%
\circ \mathcal{F}\overset{\cdot }{\rightarrow }\mathcal{K}\circ
\mathcal{G}$.

Finally, the two compositions (\ref{natCom}) and (\ref{GodCom}) of
natural transformations can be combined as
\begin{equation*}
\left( \delta \ast \gamma \right) \circ \left( \beta \ast \alpha
\right) =\left( \delta \circ \beta \right) \ast \left( \gamma
\circ \alpha \right),
\end{equation*}
where $\mathcal{A}$, $\mathcal{B}$ and $\mathcal{C}$ are three categories, $%
\mathcal{F},\mathcal{G}$, $\mathcal{H}$, $\mathcal{K}$, $\mathcal{L}$, $%
\mathcal{M}$ are six functors, and $\mathbf{\alpha }:\mathcal{F}\overset{%
\cdot }{\rightarrow }\mathcal{H}$, $\mathbf{\beta }:\mathcal{G}\overset{%
\cdot }{\rightarrow }\mathcal{K}$, $\mathbf{\gamma }:\mathcal{H}\overset{%
\cdot }{\rightarrow }\mathcal{L}$, $\mathbf{\delta }:\mathcal{K}\overset{%
\cdot }{\rightarrow }\mathcal{M}$ are four natural
transformations.

\subsection{Dinatural Transformations}
\label{dinat}

Double natural transformations are called \emph{dinatural
transformations}. An \emph{end of a functor} $S:C^{op}\times
C\rightarrow X$ is a universal dinatural transformation from a
constant $e$ to $S$. In other words, an end of $S$ is a pair
$\langle e,\omega \rangle $, where $e$ is an object of $X$ and
$\omega :e\overset{..}{\rightarrow }S$ is a \emph{wedge
(dinatural) transformation}
with the property that to every wedge $\beta :x\overset{..}{%
\rightarrow }S$ there is a unique arrow $h:x\rightarrow e$ of $B$ with $%
\beta _{c}=\omega _{c}h$ for all $a\in C$. We call $\omega $ the
\emph{ending wedge} with \emph{components} $\omega _{c}$, while
the object $e$ itself, by abuse of language, is called the end of
$S$ and written with integral notation as $\int\limits_{c}S(c,c)$;
thus
\begin{equation*}
S(c,c)\overset{\omega _{c}}{\rightarrow }\int\limits_{c}S(c,c)=e.
\end{equation*}
Note that the `variable of integration' $c$ appears twice under
the integral sign (once contravariant, once covariant) and is
`bound' by the integral
sign, in that the result no longer depends on $c$ and so is unchanged if `$c$%
' is replaced by any other letter standing for an object of the category $C$%
. These properties are like those of the letter $x$ under the usual integral symbol $%
\int f(x)\,dx$ of calculus.

Every end is manifestly a limit (see below) -- specifically, a
limit of a suitable diagram in $X$ made up of pieces like
$S(b,b)\rightarrow S(b,c)\rightarrow S(c,c)$.

For each functor $T:C\rightarrow X$ there is an isomorphism
\begin{equation*}
\int\limits_{c}S(c,c)=\int\limits_{c}Tc\cong \limfunc{Lim}T,
\end{equation*}
valid when either the end of the limit exists, carrying the ending
wedge to the limiting cone; the indicated notation thus allows us
to write any limit as an integral (an end) without explicitly
mentioning the dummy variable (the first variable $c$ of $S$).

A functor $H:X\rightarrow Y$ is said to \emph{preserve the end} of
a
functor $S:C^{op}\times C\rightarrow X$ when $\omega :e\overset{..}{%
\rightarrow }S$ an end of $S$ in $X$ implies that $H\omega :He\overset{..}{%
\rightarrow }HS$ is an and for $HS$; in symbols
\begin{equation*}
H\int\limits_{c}S(c,c)=\int\limits_{c}HS(c,c).
\end{equation*}
Similarly, $H$ \emph{creates} the end of $S$ when to each end $v:y\overset{%
..}{\rightarrow }HS$ in $Y$ there is a unique wedge $\omega :e\overset{..}{%
\rightarrow }S$ with $H\omega =v$, and this wedge $\omega $ is an
end of $S.$

The definition of the coend of a functor $S:C^{op}\times
C\rightarrow X$ is dual to that of an end. A \emph{coend} of $S$
is a pair $\langle d,\zeta
\rangle $, consisting of an object $d\in X$ and a wedge $\zeta :S\overset{..}%
{\rightarrow }d$. The object $d$ (when it exists, unique up to
isomorphism) will usually be written with an integral sign and
with the bound variable $c$ as superscript; thus
\begin{equation*}
S(c,c)\overset{\zeta _{c}}{\rightarrow }\int\limits^{c}S(c,c)=d.
\end{equation*}
The formal properties of coends are dual to those of ends. Both
are much like those for integrals in calculus (see \cite{MacLane},
for technical details).

\section{Limits and Colimits}

In abstract algebra constructions are often defined by an abstract
property which requires the existence of unique morphisms under
certain conditions. These properties are called \emph{universal
properties}. The \emph{limit} of a functor generalizes the notions
of inverse limit and product used in various parts of mathematics.
The dual notion, \emph{colimit}, generalizes direct limits and
direct sums. Limits and colimits are defined via universal
properties and provide many examples of \textit{adjoint functors}.

A \emph{limit} of a covariant functor
$\mathcal{F}:\mathcal{J}\rightarrow
\mathcal{C}$ is an object $L$ of $\mathcal{C}$, together with morphisms $%
\phi _{X}:L\rightarrow \mathcal{F}(X)$ for every object $X$ of $\mathcal{J}$%
, such that for every morphism $f:X\rightarrow Y$ in $\mathcal{J}$, we have $%
\mathcal{F}(f)\phi _{X}=\phi _{Y}$, and such that the following
\emph{universal property} is satisfied: for any object $N$ of
$\mathcal{C}$ and any set of morphisms $\psi _{X}:N\rightarrow
\mathcal{F}(X)$ such that for
every morphism $f:X\rightarrow Y$ in $\mathcal{J}$, we have $\mathcal{F}%
(f)\psi _{X}=\psi _{Y}$, there exists precisely one morphism
$u:N\rightarrow L$ such that $\phi _{X}u=\psi _{X}$ for all $X$.
If $\mathcal{F}$ has a limit (which it need not), then the limit
is defined up to a unique isomorphism, and is denoted by $\lim
\mathcal{F}$.

Analogously, a \emph{colimit} of the functor $\mathcal{F}:\mathcal{J}%
\rightarrow \mathcal{C}$ is an object $L$ of $\mathcal{C}$,
together with
morphisms $\phi _{X}:\mathcal{F}(X)\rightarrow L$ for every object $X$ of $%
\mathcal{J}$, such that for every morphism $f:X\rightarrow Y$ in
$\mathcal{J} $, we have $\phi _{Y}\mathcal{F}(X)=\phi _{X}$, and
such that the following universal property is satisfied: for any
object $N$ of $\mathcal{C}$ and any set of morphisms $\psi
_{X}:\mathcal{F}(X)\rightarrow N$ such that for every
morphism $f:X\rightarrow Y$ in $\mathcal{J}$, we have $\psi _{Y}\mathcal{F}%
(X)=\psi _{X}$, there exists precisely one morphism
$u:L\rightarrow N$ such that $u\phi _{X}=\psi _{X}$ for all $X$.
The colimit of $\mathcal{F}$,
unique up to unique isomorphism if it exists, is denoted by $\limfunc{colim}%
\mathcal{F}$.

Limits and colimits are related as follows: A functor $\mathcal{F}:\mathcal{J%
}\rightarrow \mathcal{C}$ has a colimit iff for every object $N$ of $%
\mathcal{C}$, the functor $X\longmapsto
Mor_{\mathcal{C}}(\mathcal{F}(X),N)$ (which is a covariant functor
on the dual category $\mathcal{J}^{op}$) has a
limit. If that is the case, then $Mor_{\mathcal{C}}(\limfunc{colim}\mathcal{F}%
,N)=\lim Mor_{\mathcal{C}}(\mathcal{F}(-),N)$ for every object $N$ of $%
\mathcal{C}$.

\section{Adjunction}
\label{adj}

The most important functorial operation is adjunction; as S.
MacLane once said, ``Adjoint functors arise everywhere''
\cite{MacLane}.

The \emph{adjunction} $\mathbf{\varphi }:\mathcal{F}\dashv
\mathcal{G}$ between two functors
$(\mathcal{F},\mathcal{G}):\mathcal{K}\leftrightarrows
\mathcal{L}$ of \emph{opposite variance} \cite{Kan}, represents a
\emph{weak functorial inverse:}
\begin{equation*}
\frac{f:\mathcal{F}(A)\rightarrow B}{\mathbf{\varphi
}(f):A\rightarrow \mathcal{G}(B)},
\end{equation*}
forming a \emph{natural equivalence} $\mathbf{\varphi }:\mathtt{Mor}_{%
\mathcal{K}}(\mathcal{F}(A),B) \overset{\mathbf{\varphi }}{\longrightarrow }%
\mathtt{Mor}_{\mathcal{L}}(A,\mathcal{G}(B)).$ The adjunction
isomorphism is given by a \emph{bijective correspondence} (a \emph{1--1}
and \emph{onto} map on objects) $\mathbf{\varphi
}:\mathtt{Mor}(\mathcal{K})\ni f\rightarrow \mathbf{\varphi
}(f)\in \mathtt{Mor}(\mathcal{L})$ of isomorphisms in the
two categories, $\mathcal{K}$ (with a representative object $A$), and $%
\mathcal{L}$ (with a representative object $B$). It can be
depicted as a \emph{non--commutative diagram} \\
\begin{equation*}
\bfig \putsquare<-1`1`1`1;1400`400>(-200,0)[\mathcal{F}(A)`A`B`\mathcal{G}%
(B); \mathcal{F}`f`\mathbf{\varphi }(f)`\mathcal{G}]
\put(-150,200){\oval(600,700)}\put(0,170){$\mathcal{K}$}
\put(1200,200){\oval(600,700)}\put(1000,170){$\mathcal{L}$} \efig
\end{equation*}
\bigbreak\noindent In this case $\mathcal{F}$ is called \emph{left adjoint}%
, while $\mathcal{G}$ is called \emph{right adjoint}.

In other words, an \emph{adjunction} $F\dashv G$ \emph{between two
functors} $(\mathcal{F},\mathcal{G})$ of opposite variance, from a
source
category $\mathcal{K}$ to a target category $\mathcal{L}$, is denoted by $(%
\mathcal{F},\mathcal{G},\mathbf{\eta },\mathbf{\varepsilon }):\mathcal{K}%
\leftrightarrows \mathcal{L}$. Here,
$\mathcal{F}:\mathcal{L}\rightarrow
\mathcal{K}$ is the \emph{left (upper) adjoint functor}, $%
\mathcal{G}:\mathcal{K}\leftarrow \mathcal{L}$ is the \emph{right
(lower) adjoint functor}, $\mathbf{\eta
}:1_{\mathcal{L}}\rightarrow \mathcal{G}\circ \mathcal{F}$ is the
\emph{unit natural transformation} (or, \emph{front adjunction}),
and $\mathbf{\varepsilon }:\mathcal{F}\circ \mathcal{G}\rightarrow
1_{\mathcal{K}}$ is the \emph{counit natural transformation} (or,
\emph{back adjunction}).

For example, $\mathcal{K}=\mathcal{S}$ is the category of sets and $\mathcal{%
L}=\mathcal{G}$ is the category of groups. Then $\mathcal{F}$
turns any set
into the \emph{free group} on that set, while the `forgetful'\ functor $%
\mathcal{F}^{\ast }$ turns any group into the \emph{underlying
set} of that group. Similarly, all sorts of other `free'\ and
`underlying'\ constructions are also left and right adjoints,
respectively.

Right adjoints preserve \emph{limit}\emph{s}, and left adjoints
preserve \emph{colimit}\emph{s}.

The category $\mathcal{C}$ is called a \emph{cocomplete category} if every functor $%
\mathcal{F}:\mathcal{J}\rightarrow \mathcal{C}$ has a colimit. The
following categories are cocomplete:
$\mathcal{S},\mathcal{G},\mathcal{A},\mathcal{T},$ and
$\mathcal{PT}.$

The importance of adjoint functors lies in the fact that every
functor which has a left adjoint (and therefore is a right
adjoint) is continuous. In the category $\mathcal{A}$ of Abelian
groups, this shows e.g. that the kernel of a product of
homomorphisms is naturally identified with the product of the
kernels. Also, limit functors themselves are continuous. A
covariant functor $\mathcal{F}:\mathcal{J}\rightarrow \mathcal{C}$
is \emph{co-continuous} if it transforms colimits into colimits.
Every functor which has a right adjoint (and therefore is a left adjoint) is
co-continuous.

\subsection{Application: Physiological Sensory--Motor Adjunction}

Recall that sensations from the skin,
muscles, and internal organs of the body, are transmitted to the
central nervous system via axons that enter via spinal nerves. They are called \emph{sensory pathways}. On the other hand, the motor system executes
control over the skeletal muscles of the body via several major
tracts (including pyramidal and extrapyramidal). They are called \emph{motor pathways}. Sensory--motor (or,
sensorimotor) control/coordination concerns relationships between sensation and
movement or, more broadly, between perception and action. The
interplay of sensory and motor processes provides the basis of
observable human behavior. Anatomically, its top--level, association link can be visualized
as a talk between sensory and motor Penfield's homunculi.
This sensory--motor control system can be
modelled as an adjunction between the afferent sensory functor $\mathcal{S}:%
\mathcal{BODY}\to\mathcal{BRAIN}$ and the efferent motor functor $\mathcal{M}%
:\mathcal{BRAIN}\to\mathcal{BODY}$. Thus, we have $\mathcal{SMC}%
:\mathcal{S}\dashv \mathcal{M}$, with $(\mathcal{S},\mathcal{M}):\mathcal{BRAIN}%
\leftrightarrows \mathcal{BODY}$ and depicted as

\bigbreak
\begin{equation*}
\bfig \putsquare<-1`1`1`1;1400`400>(-200,0)[\mathcal{S}(A)`A`B`\mathcal{M}%
(B); \mathcal{S}`f`\mathcal{SMC}(f)`\mathcal{M}]
\put(-150,200){\oval(700,800)}\put(-150,170){BRAIN}
\put(1200,200){\oval(850,800)}\put(850,170){BODY} \efig
\end{equation*}%
\bigbreak\bigbreak

This adjunction offers a mathematical answer to the fundamental
question: How would \emph{Nature} solve a general biodynamics
control/coordination problem? \emph{By using a weak functorial inverse of
sensory neural pathways and motor neural pathways, Nature
controls human behavior in general, and human motion in
particular}.

More generally, normal functioning of human body is achieved through interplay of a number of physiological systems -- Objects of the category BODY: musculoskeletal system, circulatory system, gastrointestinal system, integumentary system, urinary system, reproductive system, immune system and endocrine system. These systems are all interrelated, so one can say that the Morphisms between them make the proper functioning of the BODY as a whole. On the other hand, BRAIN contains the images of all above functional systems (Brain objects) and their interrelations (Brain morphisms), for the purpose of body control. This body--control performed by the brain is partly unconscious, through neuro--endocrine complex, and partly conscious, through neuro--muscular complex. A generalized sensory functor $\mathcal{SS}$ sends the information about the state of all Body objects (at any time instant) to their images in the Brain. A generalized motor functor $\mathcal{MM}$ responds to these upward sensory signals by sending downward corrective action--commands from the Brain's objects and morphisms to the Body's objects and morphisms.

For physiological details, see \cite{NatBio}. For other bio--physical applications of categorical meta-language, see \cite{GaneshSprSml,GaneshSprBig,GaneshADG}.

\section{Appendix: Groups and Related Algebraic Structures}
\label{group}

As already stated, the basic functional unit of lower biomechanics
is the special Euclidean group $SE(3)$ of rigid body motions.
In general, a \emph{group} is a pointed set $(G,e)$ with a \emph{multiplication} $%
\mu :G\times G\rightarrow G$ and an \emph{inverse} $\nu
:G\rightarrow G$ such that the following diagrams commute
\cite{Switzer}:
\begin{enumerate}
    \item \begin{equation*}
\Vtrianglepair<1`1`1`1`-1;>[G`G\times G`G`G;(e,1)`(1,e)`1`\mu `1]
\end{equation*}
($e$ is a two--sided identity)
    \item \begin{equation*}
\square <1`1`1`1;900`500>[G\times G\times G`G\times G`G\times
G`G;\mu \times 1`1\times \mu `\mu `\mu ]
\end{equation*}
(associativity)
    \item \begin{equation*}
\Vtrianglepair<1`1`1`1`-1;>[G`G\times G`G`G;(\nu ,1)`(1,\nu
)`e`\mu `e]
\end{equation*}
(inverse).
\end{enumerate}

Here $e:G\rightarrow G$ is the constant map $e(g)=e$ for all $g\in G$. $%
(e,1) $ means the map such that $(e,1)(g)=(e,g)$, etc. A group $G$
is called \emph{commutative} or \emph{Abelian group} if in
addition the following diagram commutes
\begin{equation*}
\Vtriangle<1`1`1;>[G\times G`G\times G`G;T`\mu `\mu ]
\end{equation*}
where $T:G\times G\rightarrow G\times G$ is the switch map $%
T(g_{1},g_{2})=(g_{2},g_{1}),$ for all $(g_{1},g_{2})\in G\times
G.$

A group $G$ \emph{acts}\index{group action} (on the left) on a set $A$ if there is a function $%
\alpha :G\times A\rightarrow A$ such that the following diagrams
commute \cite{Switzer}:
\begin{enumerate}
    \item \begin{equation*}
\qtriangle<1`1`1;>[A`G\times A`A;(e,1)`1`\alpha ]
\end{equation*}
    \item \begin{equation*}
\square <1`1`1`1;900`500>[G\times G\times A`G\times A`G\times
A`A;1\times \alpha `\mu \times 1`\alpha `\alpha ]
\end{equation*}%
where $(e,1)(x)=(e,x)$ for all $x\in A$. The \emph{orbit}\emph{s}
of the action are the sets $Gx=\{gx:g\in G\}$\ for all $x\in A$.
\end{enumerate}

Given two groups $(G,\ast )$ and $(H,\cdot )$, a \emph{group
homomorphism} from $(G,\ast )$ to $(H,\cdot )$ is a function
$h:G\rightarrow H$ such that for all $x$ and $y$ in $G$ it holds
that\
\begin{equation*}
h(x\ast y)=h(x)\cdot h(y).
\end{equation*}%
From this property, one can deduce that $h$ maps the identity
element $e_{G}$ of $G$ to the identity element $e_{H}$ of $H$, and
it also maps inverses to inverses in the sense that
$h(x^{-1})=h(x)^{-1}$. Hence one can say that $h$ is
\emph{compatible} with the \emph{group structure}.

The \emph{kernel} $\limfunc{Ker}h$\ of a group homomorphism $%
h:G\rightarrow H$ consists of all those elements of $G$ which are
sent by $h$ to the identity element $e_{H}$ of $H$, i.e.,
\begin{equation*}
\limfunc{Ker}h=\{x\in G:h(x)=e_{H}\}.
\end{equation*}

The \emph{image} $\limfunc{Im}h$\ of a group homomorphism
$h:G\rightarrow H$ consists of all elements of $G$ which are sent
by $h$ to $H$, i.e.,
\begin{equation*}
\limfunc{Im}h=\{h(x):x\in G\}.
\end{equation*}

The kernel is a \emph{normal subgroup} of $G$ and the image is a
\emph{subgroup} of $H$. The homomorphism $h$ is \emph{injective}
(and called a \emph{group monomorphism}) iff
$\limfunc{Ker}h=e_{G}$, i.e., iff the kernel of $h$ consists of
the identity element of $G$ only.

Similarly, a \emph{ring} (the term introduced by \emph{David Hilbert}) is a set $S$ together with two binary
operators $+$ and $\ast$ (commonly interpreted as addition and
multiplication, respectively) satisfying the following conditions:

\begin{enumerate}
\item  Additive associativity: For all $a,b,c\in S$,
$(a+b)+c=a+(b+c),$

\item  Additive commutativity: For all $a,b\in S$, $a+b=b+a,$

\item  Additive identity: There exists an element $0\in S$\ such
that for all $a\in S$, $0+a=a+0=a,$

\item  Additive inverse: For every $a\in S$\ there exists $-a\in
S$\ such that $a+(-a)=(-a)+a=0,$

\item  Multiplicative associativity: For all $a,b,c\in S$, $(a\ast
b)\ast c=a\ast (b\ast c),$

\item  Left and right distributivity: For all $a,b,c\in S$, $a\ast
(b+c)=(a\ast b)+(a\ast c)$ and $(b+c)\ast a=(b\ast a)+(c\ast a).$
\end{enumerate}

A ring is
therefore an Abelian group under addition and a semigroup under
multiplication. A ring that is commutative under multiplication,
has a unit element, and has no divisors of zero is called an
\emph{integral domain}. A ring which is also a commutative
multiplication group is called a \emph{field}. The simplest rings
are the integers $\mathbb{Z}$, polynomials $R[x]$ and $R[x,y]$ in
one and two variables, and square $n\times n$ real matrices.

An \emph{ideal} is a subset $\mathfrak{I}$ of elements in a ring
$R$ which forms an additive group and has the property that, whenever $x$ belongs to $%
R $ and $y$ belongs to $\mathfrak{I}$, then $xy$ and $yx$ belong
to $\mathfrak{I}$. For example, the set of even integers is an
ideal in the ring of integers $\mathbb{Z}$. Given an ideal
$\mathfrak{I}$, it is possible to define a factor ring
$R/\mathfrak{I}$.

A ring is called \emph{left} (respectively, \emph{right})
\emph{Noetherian} if it does not contain an infinite ascending
chain of left (respectively, right) ideals. In this case, the ring
in question is said to satisfy the ascending chain condition on
left (respectively, right) ideals. A \emph{ring} is said to be
\emph{Noetherian} if it is both left and right
Noetherian\index{Noetherian ring}. If a ring $R$ is Noetherian,
then the following are equivalent:

\begin{enumerate}
\item  $R$ satisfies the ascending chain condition on ideals.

\item  Every ideal of $R$ is finitely generated.

\item  Every set of ideals contains a maximal element.
\end{enumerate}

A \emph{module} is a mathematical object in which things can be
added together commutatively by multiplying coefficients and in
which most of the rules of manipulating vectors hold. A module is
abstractly very similar to a vector space, although in modules,
coefficients are taken in rings which are much more general
algebraic objects than the fields used in vector spaces. A
module taking its coefficients in a ring $R$ is called a module over $R$ or $%
R\mathbb{-}$module. Modules are the basic tool of homological
algebra.

Examples of modules include the set of integers $\mathbb{Z}$, the
cubic lattice in $d$ dimensions $\mathbb{Z}^{d}$, and the group ring of a group. $%
\mathbb{Z}$ is a module over itself. It is closed under addition
and subtraction. Numbers of the form $n\alpha $ for $n\in \mathbb{Z}$\ and $%
\alpha $\ a fixed integer form a submodule since, for $(n,m)\in
\mathbb{Z}$, $n\alpha \pm m\alpha =(n\pm m)\alpha $ and $(n\pm m)$
is still in $\mathbb{Z} $. Also, given two integers $a$ and $b$,
the smallest module containing $a$ and $b$ is the module for their
greatest common divisor, $\alpha =GCD(a,b)$.

A module $M$ is a \emph{Noetherian module} if it obeys the
ascending chain condition with respect to inclusion, i.e., if
every set of increasing sequences of submodules eventually becomes
constant. If a module $M$ is Noetherian, then the following are
equivalent:

\begin{enumerate}
\item  $M$ satisfies the ascending chain condition on submodules.

\item  Every submodule of $M$ is finitely generated.

\item  Every set of submodules of $M$ contains a maximal element.
\end{enumerate}

Let $I$ be a partially ordered set. A \emph{direct system} of
$R-$\emph{modules} over $I$ is an ordered pair $\{M_{i},\varphi
_{j}^{i}\}$ consisting of an indexed family of modules
$\{M_{i}:i\in I\}$ together with a family
of homomorphisms $\{\varphi _{j}^{i}:M_{i}\rightarrow M_{j}\}$ for $i\leq j$%
, such that $\varphi _{i}^{i}=1_{M_{i}}$ for all $i$ and such that
the following diagram commutes whenever $i\leq j\leq k$

\begin{equation*}
\Vtriangle<1`1`-1;>[M_{i}`M_{k}`M_{j};\varphi_{k}^{i}`
\varphi_{k}^{j}`\varphi_{j}^{i}]
\end{equation*}

Similarly, an \emph{inverse system} of $R-$\emph{modules} over $I$
is an ordered pair $\{M_{i},\psi _{i}^{j}\}$ consisting of an
indexed family of modules $\{M_{i}:i\in I\}$ together with a
family of homomorphisms $\{\psi _{i}^{j}:M_{j}\rightarrow M_{i}\}$
for $i\leq j$, such that $\psi _{i}^{i}=1_{M_{i}}$ for all $i$ and
such that the following diagram commutes whenever $i\leq j\leq k$

\begin{equation*}
\Vtriangle<1`1`-1;>[M_{k}`M_{i}`M_{j};\psi^{k}_{i}`\psi^{k}_{j}`
\psi^{j}_{i}]
\end{equation*}

\end{document}

%% file: ivdiag.tex
\newcommand{\mcm}[3]{\newcommand{#1}[#2]{{\ensuremath{#3}}}}

\usepackage{latexsym}
\usepackage{amssymb}

\mcm{\blank}{0}{(\emptybk)} \mcm{\dashbk}{0}{\mbox{---}}
\mcm{\emptybk}{0}{\:\:} \mcm{\hyph}{0}{\mbox{-}}

\mcm{\diagspace}{0}{\mbox{\hspace{2em}}}

\mcm{\cat}{1}{\mc{#1}} \mcm{\fcat}{1}{\mb{#1}}
\mcm{\mc}{1}{\mathcal{#1}} \mcm{\mr}{1}{\mathrm{#1}}
\mcm{\mi}{1}{\mathit{#1}} \mcm{\mb}{1}{\mathbf{#1}}
\mcm{\scat}{1}{\Bbb{#1}} \mcm{\twid}{1}{\widetilde{#1}}

\mcm{\elt}{0}{\in} \mcm{\sub}{0}{\,\subseteq\,}
\mcm{\such}{0}{\:|\:} \mcm{\without}{0}{\setminus}

\mcm{\atsr}{0}{\Box} \mcm{\eqv}{0}{\,\simeq\,}
\mcm{\iso}{0}{\,\cong\,}
\mcm{\of}{0}{\raisebox{0.2mm}{\ensuremath{\scriptstyle\circ}}}

\mcm{\bdry}{0}{\partial}

\mcm{\Bee}{0}{\cat{B}} \mcm{\Beep}{0}{\cat{B'}}
\mcm{\Eee}{0}{\cat{E}} \mcm{\Eeep}{0}{\cat{E'}}

\mcm{\Ess}{0}{\cat{S}} \mcm{\Tee}{0}{\cat{T}}
\mcm{\Teep}{0}{\cat{T'}} \mcm{\Stee}{0}{\scat{T}}
\mcm{\Steep}{0}{\scat{T'}}

\mcm{\blbk}{0}{\blank^{\blob}}
\mcm{\blob}{0}{\scriptscriptstyle{\bullet}}
\mcm{\stbk}{0}{\blank^{*}} \mcm{\ubl}{0}{{}^{\blob}}
\mcm{\ust}{0}{{}^{*}}

\mcm{\Cartpr}{0}{\pr{\Eee}{T}} \mcm{\Cartprp}{0}{\pr{\Eeep}{T'}}
\mcm{\Mnd}{0}{\triple{T}{\eta}{\mu}}
\mcm{\Zeropr}{0}{\pr{\Set}{\id}}

\mcm{\dopset}{0}{\ftrcat{\Delta^{\op}}{\Set}}
\mcm{\tropset}{0}{\ftrcat{\fcat{TR}^{\op}}{\Set}}

\mcm{\cod}{0}{\mr{cod}} \mcm{\dom}{0}{\mr{dom}}
\mcm{\End}{0}{\mr{End}} \mcm{\Hom}{0}{\mr{Hom}}
\mcm{\ob}{0}{\mr{ob}\,} \mcm{\op}{0}{\mr{op}}

\mcm{\comp}{0}{\mi{comp}} \mcm{\id}{0}{\mi{id}}
\mcm{\ids}{0}{\mi{ids}} \mcm{\mult}{0}{\mi{mult}}
\mcm{\unit}{0}{\mi{unit}}

\mcm{\Ab}{0}{\fcat{Ab}} \mcm{\Alg}{0}{\fcat{Alg}}
\mcm{\Bim}{1}{\fcat{Bim}(#1)} \mcm{\Cat}{0}{\fcat{Cat}}
\mcm{\Cay}{0}{\fcat{Cay}} \mcm{\Cpn}{1}{\pr{\Set/S_{#1}}{T_{#1}}}
\mcm{\fc}{0}{\fcat{fc}} \mcm{\fm}{0}{\fcat{fm}}
\mcm{\Graph}{0}{\fcat{Graph}} \mcm{\Gy}{0}{\fcat{Gy}}
\mcm{\Hpn}{1}{\pr{\Eee_{#1}}{P_{#1}}} \mcm{\Mon}{0}{\mb{Mon}}
\mcm{\Multicat}{0}{\fcat{Multicat}} \mcm{\One}{0}{\fcat{1}}
\mcm{\PD}{1}{\fcat{PD}_{#1}} \mcm{\Prof}{0}{\fcat{Prof}}
\mcm{\Set}{0}{\fcat{Set}} \mcm{\Span}{0}{\fcat{Span}}
\mcm{\Ssq}{0}{\fcat{Ssq}} \mcm{\Struc}{0}{\fcat{Struc}}
\mcm{\Sym}{0}{\fcat{Sym}} \mcm{\TR}{1}{\fcat{TR}(#1)}
\mcm{\Tr}{0}{\fcat{Tr}} \mcm{\Twocat}{0}{\fcat{2\hyph\Cat}}

\mcm{\integers}{0}{\mathbb{Z}}

\mcm{\range}{2}{#1,\,\ldots\,,#2}
\mcm{\bftuple}{2}{\tuplebts{\range{#1}{#2}}}
\mcm{\tuple}{3}{\tuplebts{\range{#1,#2}{#3}}}
\mcm{\rttuple}{1}{\tuplebts{\,\ldots\,,#1}}
\mcm{\abftuple}{2}{\atuplebts{\range{#1}{#2}}}
\mcm{\atuple}{3}{\atuplebts{\range{#1,#2}{#3}}}
\mcm{\arttuple}{1}{\atuplebts{\,\ldots\,,#1}}
\mcm{\sqbftuple}{2}{\obt\range{#1}{#2}\cbt}
\mcm{\pr}{2}{\tuplebts{#1,#2}}
\mcm{\triple}{3}{\tuplebts{#1,#2,#3}}

\mcm{\eend}{2}{#1[#2]} \mcm{\ehom}{3}{#1[#2,#3]}
\mcm{\ftrcat}{2}{[#1,#2]} \mcm{\homset}{3}{#1(#2,#3)}
\mcm{\multihom}{3}{#1(#2;#3)}
\mcm{\relhom}{5}{#1_{#2}(\range{#3}{#4};#5)}

\mcm{\go}{0}{\rTo} \mcm{\goby}{1}{\rTo^{#1}}
\mcm{\goesto}{0}{\,\longmapsto\,} \mcm{\goiso}{0}{\goby{\diso}}
\mcm{\monic}{0}{\rMonic} \mcm{\og}{0}{\lTo}
\mcm{\ogby}{1}{\lTo^{#1}}

\mcm{\gph}{2}{\spn{#1}{T #2}{#2}} \mcm{\graph}{4}{\spaan{#1}{T
#2}{#2}{#3}{#4}} \mcm{\oppair}{2}{\stackrel{\rTo^{#1}}{\lTo_{#2}}}
\mcm{\parpair}{2}{\stackrel{\rTo^{#1}}{\rTo_{#2}}}
\mcm{\spn}{3}{#2 \og #1 \go #3} \mcm{\spaan}{5}{#2 \ogby{#4} #1
\goby{#5} #3}

\mcm{\bktdvslob}{3}
    {\left(
    \begin{diagram}[height=1.5em]
    #1      \\
    \dTo>{\,#2} \\
    #3      \\
    \end{diagram}
    \right)}
\mcm{\slob}{3}{(#1 \goby{#2} #3)} \mcm{\vslob}{3}
    {\left.
    \begin{diagram}[height=1.5em]
    #1      \\
    \dTo>{\,#2} \\
    #3      \\
    \end{diagram}
    \right.}

\newenvironment{tree}
    {\begin{diagram}[height=1em,width=.75em,abut,noPS,tight]}
    {\end{diagram}}

\mcm{\enode}{0}{\circ}

\mcm{\nl}{1}{\stackrel{\textstyle #1}{\node}}
\mcm{\node}{0}{\bullet}

\mcm{\utree}{0}{\node}

\mcm{\diso}{0}{\sim}

\mcm{\vdiso}{0}{\wr}

\mcm{\nat}{0}{\mathbb{N}}

\mcm{\Onepr}{0}{\pr{\Graph}{\fc}}
\newlength{\nllwidth}
\newlength{\nllheight}
\newcommand{\stackbelow}[2]{%
\settowidth{\nllwidth}{\ensuremath{#1}\ensuremath{#2}}%
\settoheight{\nllheight}{\ensuremath{#2}}%
\addtolength{\nllheight}{.3ex}%
\mbox{%
\ensuremath{#1}%
\hspace{-.5\nllwidth}%
\raisebox{-1\nllheight}{\ensuremath{#2}}}}
\mcm{\nlal}{2}{\stackbelow{\nl{#1}}{#2}}
\mcm{\nll}{1}{\stackbelow{\node}{#1}} \mcm{\wun}{0}{\fcat{1}}
\mcm{\atuplebts}{1}{\langle #1 \rangle} \mcm{\tuplebts}{1}{(#1)}
\mcm{\bo}{0}{(} \mcm{\bc}{0}{)}
\mcm{\UBilax}{0}{\fcat{UBicat}_\mr{lax}}
\mcm{\UBiwk}{0}{\fcat{UBicat}_\mr{wk}}
\mcm{\UBistr}{0}{\fcat{UBicat}_\mr{str}}
\mcm{\Bilax}{0}{\fcat{Bicat}_\mr{lax}}
\mcm{\Biwk}{0}{\fcat{Bicat}_\mr{wk}}
\mcm{\Bistr}{0}{\fcat{Bicat}_\mr{str}} \mcm{\rotsub}{0}{\cup
\raisebox{0.1em}{$\scriptstyle{|}$}} \mcm{\pd}{0}{\fcat{pd}}
\mcm{\rep}{1}{\widehat{#1}} \mcm{\ovln}{1}{\overline{#1}}
\mcm{\Gph}{0}{\fcat{Gph}} \mcm{\tr}{0}{\fcat{tr}}

\mcm{\ladj}{0}{\,\dashv\,} \mcm{\zeropd}{0}{\node}
    {\end{diagram}}
\mcm{\END}{0}{\fcat{End}} \mcm{\HOM}{0}{\fcat{Hom}}

\newlength{\gwidth} 
\newlength{\gvert}  
\newlength{\gdrop}  
\newlength{\gbaredrop}  
\newlength{\goffset}    
\newlength{\gtemp}  

\newcommand{\present}[1]{%
\makebox[1\gwidth]{%
\rule[-1\gdrop]{0ex}{1\gvert}%
\raisebox{-1\gbaredrop}{#1}}}

\newcommand{\presentl}[1]{%
\makebox[1\gwidth][l]{%
\rule[-1\gdrop]{0ex}{1\gvert}%
\raisebox{-1\gbaredrop}{#1}}}

\newcommand{\presentr}[1]{%
\makebox[1\gwidth][r]{%
\rule[-1\gdrop]{0ex}{1\gvert}%
\raisebox{-1\gbaredrop}{#1}}}

\newcommand{\ginitdims}[2]{
\setlength{\unitlength}{1em}
\setlength{\goffset}{.25\unitlength}
\setlength{\gwidth}{#1\unitlength}
\setlength{\gvert}{#2\unitlength}
\setlength{\gdrop}{.5\gvert}
\addtolength{\gdrop}{-1\goffset}
\setlength{\gbaredrop}{1\gdrop}
\addtolength{\gvert}{.6\unitlength}
\addtolength{\gdrop}{.3\unitlength}}    

\newcommand{\cinitdims}[2]{
\setlength{\unitlength}{1em}
\setlength{\goffset}{.35\unitlength}
\setlength{\gwidth}{#1\unitlength}
\setlength{\gvert}{#2\unitlength}
\setlength{\gdrop}{.5\gvert}
\addtolength{\gdrop}{-1\goffset}
\setlength{\gbaredrop}{1\gdrop}
\addtolength{\gvert}{.6\unitlength}
\addtolength{\gdrop}{.3\unitlength}}    

\newcommand{\gsinitdims}[2]{
\setlength{\unitlength}{0.5em}
\setlength{\goffset}{.25\unitlength}
\setlength{\gwidth}{#1\unitlength}
\setlength{\gvert}{#2\unitlength}
\setlength{\gdrop}{.5\gvert}
\addtolength{\gdrop}{-1\goffset}
\setlength{\gbaredrop}{1\gdrop}
\addtolength{\gvert}{.6\unitlength}
\addtolength{\gdrop}{.3\unitlength}}    

\newcommand{\sidespic}[1]{%
\settowidth{\gtemp}{\ensuremath{#1}}%
\addtolength{\gwidth}{1\gtemp}}

\newcommand{\abovepic}[1]{%
\settoheight{\gtemp}{\ensuremath{#1}}%
\addtolength{\gvert}{1\gtemp}%
\settodepth{\gtemp}{\ensuremath{#1}}%
\addtolength{\gvert}{1\gtemp}}

\newcommand{\belowpic}[1]{%
\settoheight{\gtemp}{\ensuremath{#1}}%
\addtolength{\gvert}{1\gtemp}%
\addtolength{\gdrop}{1\gtemp}%
\settodepth{\gtemp}{\ensuremath{#1}}%
\addtolength{\gvert}{1\gtemp}%
\addtolength{\gdrop}{1\gtemp}}

\newcommand{\cell}[4]{\put(#1,#2){\makebox(0,0)[#3]{\ensuremath{#4}}}}
\mcm{\zmark}{0}{\scriptstyle{\bullet}}

\newcommand{\pregfst}[1]{%
\begin{picture}(0.5,0.2)(-0.5,-0.2)%
\cell{-0.1}{-0.2}{tr}{#1}%
\cell{0}{0}{c}{\zmark}%
\end{picture}}

\mcm{\gfst}{1}{%
\ginitdims{0.5}{0.4}%
\sidespic{#1}%
\belowpic{#1}%
\presentr{\pregfst{#1}}}

\newcommand{\preglst}[1]{%
\begin{picture}(0.5,0.2)(0,-0.2)%
\cell{0.1}{-0.2}{tl}{#1}%
\cell{0.05}{0}{c}{\zmark}%
\end{picture}}

\mcm{\glst}{1}{%
\ginitdims{.5}{.4}%
\sidespic{#1}%
\belowpic{#1}%
\presentl{\preglst{#1}}}

\newcommand{\preglft}[1]{%
\begin{picture}(0,0.2)(0,-0.2)%
\cell{-0.1}{-0.2}{tr}{#1}%
\cell{0.05}{0}{c}{\zmark}%
\end{picture}}

\mcm{\glft}{1}{%
\ginitdims{0}{.4}%
\belowpic{#1}%
\present{\preglft{#1}}}

\newcommand{\pregrgt}[1]{%
\begin{picture}(0,0.2)(0,-0.2)%
\cell{0.1}{-0.2}{tl}{#1}%
\cell{0.05}{0}{c}{\zmark}%
\end{picture}}

\mcm{\grgt}{1}{%
\ginitdims{0}{.4}%
\belowpic{#1}%
\present{\pregrgt{#1}}}

\newcommand{\pregblw}[1]{%
\begin{picture}(0,0.3)(0,-0.3)
\cell{0}{-0.3}{t}{#1}%
\cell{0.05}{0}{c}{\zmark}%
\end{picture}}

\mcm{\gblw}{1}{%
\ginitdims{0}{.6}%
\belowpic{#1}%
\present{\pregblw{#1}}}

\newcommand{\pregfbw}[1]{%
\begin{picture}(0,0.65)(0,-0.65)
\cell{0}{-0.65}{t}{#1}%
\cell{0.05}{0}{c}{\zmark}%
\end{picture}}

\mcm{\gfbw}{1}{%
\ginitdims{0}{1.3}%
\belowpic{#1}%
\present{\pregfbw{#1}}}

\newcommand{\pregzero}[1]{%
\begin{picture}(0.8,0.4)(-0.4,-0.4)
\cell{0}{-0.4}{t}{#1}%
\cell{0}{0}{c}{\zmark}%
\end{picture}}

\mcm{\gzero}{1}{%
\ginitdims{0.8}{.6}%
\belowpic{#1}%
\sidespic{#1}%
\present{\pregzero{#1}}}

\newcommand{\pregone}[1]{%
\begin{picture}(5,0.4)(0,-0.2)%
\cell{2.5}{0.2}{b}{#1}%
\put(0,0){\vector(1,0){5}}%
\end{picture}}

\mcm{\gone}{1}{%
\ginitdims{5}{0.4}%
\abovepic{#1}%
\present{\pregone{#1}}}

\newcommand{\pregtwo}[3]{%
\begin{picture}(5,3.4)(0,-0.2)%
\cell{2.5}{3.2}{b}{#1}%
\cell{2.5}{-.2}{t}{#2}%
\cell{2.7}{1.5}{l}{#3}%
\qbezier(0,1.5)(2.5,4.5)(5,1.5)%
\qbezier(0,1.5)(2.5,-1.5)(5,1.5)%
\put(5,1.5){\vector(1,-1){0}}%
\put(5,1.5){\vector(1,1){0}}%
\put(2.5,2.5){\vector(0,-1){2}}%
\end{picture}}

\mcm{\gtwo}{3}{%
\ginitdims{5}{3.4}%
\abovepic{#1}%
\belowpic{#2}%
\present{\pregtwo{#1}{#2}{#3}}}

\newcommand{\pregthree}[5]{%
\begin{picture}(5,5.4)(0,-1.2)%
\cell{2.5}{4.2}{b}{#1}%
\cell{1.5}{1.7}{b}{#2}%
\cell{2.5}{-1.2}{t}{#3}%
\cell{2.7}{2.75}{l}{#4}%
\cell{2.7}{0.25}{l}{#5}%
\qbezier(0,1.5)(2.5,6.5)(5,1.5)%
\qbezier(0,1.5)(2.5,-3.5)(5,1.5)%
\put(0,1.5){\vector(1,0){5}}%
\put(2.5,3.5){\vector(0,-1){1.5}}%
\put(2.5,1){\vector(0,-1){1.5}}%
\put(5,1.5){\vector(1,-3){0}}%
\put(5,1.5){\vector(1,3){0}}%
\end{picture}}

\mcm{\gthree}{5}{%
\ginitdims{5}{5.4}%
\abovepic{#1}%
\belowpic{#3}%
\present{\pregthree{#1}{#2}{#3}{#4}{#5}}}

\newcommand{\pregfour}[7]{%
\begin{picture}(5,8.4)(0,-2.7)%
\cell{2.5}{5.7}{b}{#1}%
\cell{1.5}{2.8}{b}{#2}%
\cell{1.5}{0.2}{t}{#3}%
\cell{2.5}{-2.7}{t}{#4}%
\cell{2.7}{4.25}{l}{#5}%
\cell{2.7}{1.5}{l}{#6}%
\cell{2.7}{-1.25}{l}{#7}%
\qbezier(0,1.5)(2.5,9.5)(5,1.5)%
\qbezier(0,1.5)(2.5,4)(5,1.5)%
\qbezier(0,1.5)(2.5,-1)(5,1.5)%
\qbezier(0,1.5)(2.5,-6.5)(5,1.5)%
\put(2.5,5.25){\vector(0,-1){2}}%
\put(2.5,2.5){\vector(0,-1){2}}%
\put(2.5,-0.25){\vector(0,-1){2}}%
\put(5,1.5){\vector(1,-4){0}}%
\put(5,1.5){\vector(4,-3){0}}%
\put(5,1.5){\vector(4,3){0}}%
\put(5,1.5){\vector(1,4){0}}%
\end{picture}}

\mcm{\gfour}{7}{%
\ginitdims{5}{8.4}%
\abovepic{#1}%
\belowpic{#4}%
\present{\pregfour{#1}{#2}{#3}{#4}{#5}{#6}{#7}}}

\newcommand{\pregthreecell}[5]{%
\begin{picture}(8,5)(-4,-2.5)%
\cell{0}{2.5}{b}{#1}%
\cell{0}{-2.5}{t}{#2}%
\cell{-1.7}{0}{r}{#3}%
\cell{1.7}{0}{l}{#4}%
\cell{0}{0.2}{b}{#5}%
\qbezier(-4,0)(0,4.2)(4,0)%
\qbezier(-4,0)(0,-4.2)(4,0)%
\qbezier(-0.5,1.8)(-2.5,0)(-0.5,-1.8)%
\qbezier(0.5,1.8)(2.5,0)(0.5,-1.8)%
\put(-1,0){\vector(1,0){2}}%
\put(4,0){\vector(1,-1){0}}%
\put(4,0){\vector(1,1){0}}%
\put(-0.5,-1.8){\vector(1,-1){0}}%
\put(0.5,-1.8){\vector(-1,-1){0}}%
\end{picture}}

\mcm{\gthreecell}{5}{%
\ginitdims{8}{5}%
\abovepic{#1}%
\belowpic{#2}%
\present{\pregthreecell{#1}{#2}{#3}{#4}{#5}}}

\newcommand{\pregthreecellu}{%
\begin{picture}(5,3.4)(-0.5,-0.2)%
\qbezier(-.5,1.5)(2,4.5)(4.5,1.5)%
\qbezier(-.5,1.5)(2,-1.5)(4.5,1.5)%
\qbezier(1.5,2.7)(0.5,1.5)(1.5,0.3)%
\qbezier(2.5,2.7)(3.5,1.5)(2.5,0.3)%
\put(1.3,1.5){\vector(1,0){1.4}}%
\put(4.5,1.5){\vector(1,-1){0}}%
\put(4.5,1.5){\vector(1,1){0}}%
\put(1.5,0.3){\vector(2,-3){0}}%
\put(2.5,0.3){\vector(-2,-3){0}}%
\end{picture}}

\mcm{\gthreecellu}{0}{%
\ginitdims{5}{3.4}%
\present{\pregthreecellu}}

\newcommand{\pregtwocentre}[3]{%
\begin{picture}(5,3.4)(0,-0.2)%
\cell{2.5}{3.2}{b}{#1}%
\cell{2.5}{-.2}{t}{#2}%
\cell{2.5}{1.5}{c}{#3}%
\qbezier(0,1.5)(2.5,4.5)(5,1.5)%
\qbezier(0,1.5)(2.5,-1.5)(5,1.5)%
\put(5,1.5){\vector(1,-1){0}}%
\put(5,1.5){\vector(1,1){0}}%
\put(2.5,2.5){\vector(0,-1){2}}%
\end{picture}}

\mcm{\gtwocentre}{3}{%
\ginitdims{5}{3.4}%
\abovepic{#1}%
\belowpic{#2}%
\present{\pregtwocentre{#1}{#2}{#3}}}

\newcommand{\pregspecialone}[9]{%
\begin{picture}(8,8)(-4,-4)%
\cell{0}{3.9}{b}{#1}%
\cell{-2}{-0.2}{t}{#2}%
\cell{0}{-3.9}{t}{#3}%
\cell{-1.5}{1.1}{r}{#4}%
\cell{0.2}{1.5}{l}{#5}%
\cell{1.5}{1.1}{l}{#6}%
\cell{0.2}{-2}{l}{#7}%
\cell{-0.9}{2.3}{b}{#8}%
\cell{0.9}{2.3}{b}{#9}%
\qbezier(-4,0)(0,8)(4,0)%
\qbezier(-4,0)(0,-8)(4,0)%
\qbezier(-0.5,3.4)(-3.5,2)(-0.5,0.6)%
\qbezier(0.5,3.4)(3.5,2)(0.5,0.6)%
\put(-4,0){\vector(1,0){8}}%
\put(0,3.4){\vector(0,-1){2.8}}%
\put(0,-0.8){\vector(0,-1){2.4}}%
\put(-1.5,2.2){\vector(1,0){1.2}}%
\put(0.3,2.2){\vector(1,0){1.2}}%
\put(4,0){\vector(1,-2){0}}%
\put(4,0){\vector(1,2){0}}%
\put(-0.5,0.6){\vector(2,-1){0}}%
\put(0.5,0.6){\vector(-2,-1){0}}%
\end{picture}}

\mcm{\gspecialone}{9}{%
\ginitdims{8}{8}%
\abovepic{#1}%
\belowpic{#3}%
\present{\pregspecialone{#1}{#2}{#3}{#4}{#5}{#6}{#7}{#8}{#9}}}

\newcommand{\pregspecialtwo}{%
\begin{picture}(5,3.4)(0,-0.2)%
\qbezier(0,1.5)(2.5,4.5)(5,1.5)%
\qbezier(0,1.5)(2.5,-1.5)(5,1.5)%
\qbezier(1.7,2.5)(0,1.5)(1.7,0.5)%
\qbezier(3.3,2.5)(5,1.5)(3.3,0.5)%
\put(5,1.5){\vector(1,-1){0}}%
\put(5,1.5){\vector(1,1){0}}%
\put(1.7,0.5){\vector(3,-2){0}}%
\put(3.3,0.5){\vector(-3,-2){0}}%
\put(2.5,2.5){\vector(0,-1){2}}%
\put(1.2,1.5){\vector(1,0){1}}%
\put(2.8,1.5){\vector(1,0){1}}%
\end{picture}}

\mcm{\gspecialtwo}{0}{%
\ginitdims{5}{3.4}%
\present{\pregspecialtwo}}

\newcommand{\pregspecialthree}{%
\begin{picture}(5,5.4)(0,-1.2)%
\qbezier(0,1.5)(2.5,6.5)(5,1.5)%
\qbezier(0,1.5)(2.5,-3.5)(5,1.5)%
\qbezier(2,3.5)(1,2.75)(2,2)%
\qbezier(3,3.5)(4,2.75)(3,2)%
\qbezier(2,1)(1,0.25)(2,-0.5)%
\qbezier(3,1)(4,0.25)(3,-0.5)%
\put(0,1.5){\vector(1,0){5}}%
\put(1.5,2.75){\vector(1,0){2}}%
\put(1.5,0.25){\vector(1,0){2}}%
\put(5,1.5){\vector(1,-3){0}}%
\put(5,1.5){\vector(1,3){0}}%
\put(2,2){\vector(1,-1){0}}%
\put(3,2){\vector(-1,-1){0}}%
\put(2,-0.5){\vector(1,-1){0}}%
\put(3,-0.5){\vector(-1,-1){0}}%
\end{picture}}

\mcm{\gspecialthree}{0}{%
\ginitdims{5}{5.4}%
\present{\pregspecialthree}}

\newcommand{\pregonew}[1]{%
\begin{picture}(8,0.4)(0,-0.2)%
\cell{4}{0.2}{b}{#1}%
\put(0,0){\vector(1,0){8}}%
\end{picture}}

\mcm{\gonew}{1}{%
\ginitdims{8}{0.4}%
\abovepic{#1}%
\present{\pregonew{#1}}}

\mcm{\gzersu}{0}{%
\gsinitdims{0}{.6}%
\present{\pregblw{}}}

\mcm{\gonesu}{0}{%
\gsinitdims{5}{0.4}%
\present{\pregone{}}}

\mcm{\gtwosu}{0}{%
\gsinitdims{5}{3.4}%
\present{\pregtwo{}{}{}}}

\mcm{\gthreesu}{0}{%
\gsinitdims{5}{5.4}%
\present{\pregthree{}{}{}{}{}}}

\mcm{\gfoursu}{0}{%
\gsinitdims{5}{8.4}%
\present{\pregfour{}{}{}{}{}{}{}}}

\newcommand{\precone}[1]{%
\begin{picture}(4.2,0.4)(-0.3,-0.2)%
\cell{1.8}{0.2}{b}{#1}%
\put(0,0){\vector(1,0){3.6}}%
\end{picture}}

\mcm{\cone}{1}{%
\cinitdims{4.2}{0.4}%
\abovepic{#1}%
\present{\precone{#1}}}

\mcm{\gfstsu}{0}{%
\gsinitdims{0.5}{0.4}%
\presentr{\pregfst{}}}

\mcm{\glstsu}{0}{%
\gsinitdims{0.5}{0.4}%
\presentl{\preglst{}}}

\newcommand{\prectwodbl}[3]%
{\begin{picture}(4.2,3.4)(-0.1,-0.2)%
\cell{2}{3.2}{b}{#1}%
\cell{2}{-0.2}{t}{#2}%
\cell{2.3}{1.5}{l}{#3}%
\qbezier(0,2)(2,4)(4,2)%
\qbezier(0,1)(2,-1)(4,1)%
\put(4,2){\vector(1,-1){0}}%
\put(4,1){\vector(1,1){0}}%
\put(1.9,2.5){\line(0,-1){1.8}}%
\put(2.1,2.5){\line(0,-1){1.8}}%
\cell{2.01}{0.4}{b}{\vee}%
\end{picture}}

\mcm{\ctwodbl}{3}{%
\cinitdims{4.2}{3.4}%
\abovepic{#1}%
\belowpic{#2}%
\present{\prectwodbl{#1}{#2}{#3}}}

\newcommand{\precthreedbl}[5]{%
\begin{picture}(4.2,5.4)(-0.1,-0.2)%
\cell{2}{5.2}{b}{#1}%
\cell{1}{2.7}{b}{#2}%
\cell{2}{-.2}{t}{#3}%
\cell{2.3}{3.75}{l}{#4}%
\cell{2.3}{1.25}{l}{#5}%
\qbezier(0,3)(2,7)(4,3)%
\qbezier(0,2)(2,-2)(4,2)%
\put(0,2.5){\vector(1,0){4}}%
\put(1.9,4.5){\line(0,-1){1.3}}%
\put(2.1,4.5){\line(0,-1){1.3}}%
\cell{2.01}{2.9}{b}{\vee}%
\put(1.9,2){\line(0,-1){1.3}}%
\put(2.1,2){\line(0,-1){1.3}}%
\cell{2.01}{0.4}{b}{\vee}%
\put(4,3){\vector(1,-3){0}}%
\put(4,2){\vector(1,3){0}}%
\end{picture}}

\mcm{\cthreedbl}{5}{%
\cinitdims{4.2}{5.4}%
\abovepic{#1}%
\belowpic{#3}%
\present{\precthreedbl{#1}{#2}{#3}{#4}{#5}}}

\newcommand{\precthreecelltrp}[5]{%
\begin{picture}(8.2,5)(-4.1,-2.5)%
\cell{0}{2.5}{b}{#1}%
\cell{0}{-2.5}{t}{#2}%
\cell{-1.8}{0}{r}{#3}%
\cell{1.8}{0}{l}{#4}%
\cell{0}{0.3}{b}{#5}%
\qbezier(-4,0.5)(0,4)(4,0.5)%
\qbezier(-4,-0.5)(0,-4)(4,-0.5)%
\qbezier(-0.6,2)(-2.6,0)(-0.6,-2)%
\qbezier(-0.4,2)(-2.4,0)(-0.5,-1.9)%
\cell{-0.6}{-2}{b}{\lrcorner}%
\qbezier(0.4,2)(2.4,0)(0.5,-1.9)%
\qbezier(0.6,2)(2.6,0)(0.6,-2)%
\cell{0.65}{-2}{b}{\llcorner}%
\put(-1,0.15){\line(1,0){1.7}}%
\put(-1,0){\line(1,0){2}}%
\put(-1,-0.15){\line(1,0){1.7}}%
\cell{1.15}{0}{r}{>}%
\put(4,0.5){\vector(1,-1){0}}%
\put(4,-0.5){\vector(1,1){0}}%
\end{picture}}

\mcm{\cthreecelltrp}{5}{%
\cinitdims{8.2}{5}%
\abovepic{#1}%
\belowpic{#2}%
\present{\precthreecelltrp{#1}{#2}{#3}{#4}{#5}}}

\newcommand{\prectwo}[3]%
{\begin{picture}(4.2,3.4)(-0.1,-0.2)%
\cell{2}{3.2}{b}{#1}%
\cell{2}{-0.2}{t}{#2}%
\cell{2.2}{1.5}{l}{#3}%
\qbezier(0,2)(2,4)(4,2)%
\qbezier(0,1)(2,-1)(4,1)%
\put(4,2){\vector(1,-1){0}}%
\put(4,1){\vector(1,1){0}}%
\put(2,2.5){\vector(0,-1){2}}%
\end{picture}}

\mcm{\ctwo}{3}{%
\cinitdims{4.2}{3.4}%
\abovepic{#1}%
\belowpic{#2}%
\present{\prectwo{#1}{#2}{#3}}}

\newcommand{\precthree}[5]{%
\begin{picture}(4.2,5.4)(-0.1,-0.2)%
\cell{2}{5.2}{b}{#1}%
\cell{1}{2.7}{b}{#2}%
\cell{2}{-.2}{t}{#3}%
\cell{2.2}{3.75}{l}{#4}%
\cell{2.2}{1.25}{l}{#5}%
\qbezier(0,3)(2,7)(4,3)%
\qbezier(0,2)(2,-2)(4,2)%
\put(0,2.5){\vector(1,0){4}}%
\put(2,4.5){\vector(0,-1){1.5}}%
\put(2,2){\vector(0,-1){1.5}}%
\put(4,3){\vector(1,-3){0}}%
\put(4,2){\vector(1,3){0}}%
\end{picture}}

\mcm{\cthree}{5}{%
\cinitdims{4.2}{5.4}%
\abovepic{#1}%
\belowpic{#3}%
\present{\precthree{#1}{#2}{#3}{#4}{#5}}}

\newcommand{\prectwoop}[3]%
{\begin{picture}(4.2,3.4)(-0.1,-0.2)%
\cell{2}{3.2}{b}{#1}%
\cell{2}{-0.2}{t}{#2}%
\cell{2.2}{1.5}{l}{#3}%
\qbezier(0,2)(2,4)(4,2)%
\qbezier(0,1)(2,-1)(4,1)%
\put(0,2){\vector(-1,-1){0}}%
\put(0,1){\vector(-1,1){0}}%
\put(2,2.5){\vector(0,-1){2}}%
\end{picture}}

\mcm{\ctwoop}{3}{%
\cinitdims{4.2}{3.4}%
\abovepic{#1}%
\belowpic{#2}%
\present{\prectwoop{#1}{#2}{#3}}}

\newcommand{\prectwopar}[4]{%
\begin{picture}(4.2,3.4)(-0.1,-0.2)%
\cell{2}{3.2}{b}{#1}%
\cell{2}{-0.2}{t}{#2}%
\cell{1.6}{1.5}{r}{#3}%
\cell{2.4}{1.5}{l}{#4}%
\qbezier(0,2)(2,4)(4,2)%
\qbezier(0,1)(2,-1)(4,1)%
\put(4,2){\vector(1,-1){0}}%
\put(4,1){\vector(1,1){0}}%
\put(1.8,2.5){\vector(0,-1){2}}%
\put(2.2,2.5){\vector(0,-1){2}}%
\end{picture}}

\mcm{\ctwopar}{4}{%
\cinitdims{4.2}{3.4}%
\abovepic{#1}%
\belowpic{#2}%
\present{\prectwopar{#1}{#2}{#3}{#4}}}

\newcommand{\precthreein}[5]{%
\begin{picture}(4.2,5.4)(-0.1,-0.2)%
\cell{2}{5.2}{b}{#1}%
\cell{1}{2.7}{b}{#2}%
\cell{2}{-.2}{t}{#3}%
\cell{2.2}{3.75}{l}{#4}%
\cell{2.2}{1.25}{l}{#5}%
\qbezier(0,3)(2,7)(4,3)%
\qbezier(0,2)(2,-2)(4,2)%
\put(0,2.5){\vector(1,0){4}}%
\put(2,4.5){\vector(0,-1){1.5}}%
\put(2,0.5){\vector(0,1){1.5}}%
\put(4,3){\vector(1,-3){0}}%
\put(4,2){\vector(1,3){0}}%
\end{picture}}

\mcm{\cthreein}{5}{%
\cinitdims{4.2}{5.4}%
\abovepic{#1}%
\belowpic{#3}%
\present{\precthreein{#1}{#2}{#3}{#4}{#5}}}

\newcommand{\precthreecell}[5]{%
\begin{picture}(8.2,5)(-4.1,-2.5)%
\cell{0}{2.5}{b}{#1}%
\cell{0}{-2.5}{t}{#2}%
\cell{-1.7}{0}{r}{#3}%
\cell{1.7}{0}{l}{#4}%
\cell{0}{0.2}{b}{#5}%
\qbezier(-4,0.5)(0,4)(4,0.5)%
\qbezier(-4,-0.5)(0,-4)(4,-0.5)%
\qbezier(-0.5,2)(-2.5,0)(-0.5,-2)%
\qbezier(0.5,2)(2.5,0)(0.5,-2)%
\put(-1,0){\vector(1,0){2}}%
\put(4,0.5){\vector(1,-1){0}}%
\put(4,-0.5){\vector(1,1){0}}%
\put(-0.5,-2){\vector(1,-1){0}}%
\put(0.5,-2){\vector(-1,-1){0}}%
\end{picture}}

\mcm{\cthreecell}{5}{%
\cinitdims{8.2}{5}%
\abovepic{#1}%
\belowpic{#2}%
\present{\precthreecell{#1}{#2}{#3}{#4}{#5}}}

\newcommand{\precthreecellpar}[6]{%
\begin{picture}(8.2,5)(-4.1,-2.5)%
\cell{0}{2.5}{b}{#1}%
\cell{0}{-2.5}{t}{#2}%
\cell{-1.7}{0}{r}{#3}%
\cell{1.7}{0}{l}{#4}%
\cell{0}{0.4}{b}{#5}%
\cell{0}{-0.4}{t}{#6}%
\qbezier(-4,0.5)(0,4)(4,0.5)%
\qbezier(-4,-0.5)(0,-4)(4,-0.5)%
\qbezier(-0.5,2)(-2.5,0)(-0.5,-2)%
\qbezier(0.5,2)(2.5,0)(0.5,-2)%
\put(-1,0.2){\vector(1,0){2}}%
\put(-1,-0.2){\vector(1,0){2}}%
\put(4,0.5){\vector(1,-1){0}}%
\put(4,-0.5){\vector(1,1){0}}%
\put(-0.5,-2){\vector(1,-1){0}}%
\put(0.5,-2){\vector(-1,-1){0}}%
\end{picture}}

\mcm{\cthreecellpar}{6}{%
\cinitdims{8.2}{5}%
\abovepic{#1}%
\belowpic{#2}%
\present{\precthreecellpar{#1}{#2}{#3}{#4}{#5}{#6}}}

\newcommand{\prectwov}[5]{%
\begin{picture}(3.4,4.2)(0.8,0.9)%
\cell{2.5}{5.1}{b}{#1}%
\cell{2.5}{0.9}{t}{#2}%
\cell{0.8}{3}{r}{#3}%
\cell{4.2}{3}{l}{#4}%
\cell{2.5}{3.2}{b}{#5}%
\qbezier(2,5)(0,3)(2,1)%
\qbezier(3,5)(5,3)(3,1)%
\put(2,1){\vector(1,-1){0}}%
\put(3,1){\vector(-1,-1){0}}%
\put(1.5,3){\vector(1,0){2}}%
\end{picture}}

\mcm{\ctwov}{5}{%
\cinitdims{3.4}{4.2}%
\abovepic{#1}%
\belowpic{#2}%
\sidespic{#3}%
\sidespic{#4}%
\present{\prectwov{#1}{#2}{#3}{#4}{#5}}}

\newcommand{\precthreecellv}[7]{%
\begin{picture}(5,8.2)(0.5,-1.6)%
\cell{3}{6.6}{b}{#1}%
\cell{3}{-1.6}{t}{#2}%
\cell{0.5}{2.5}{r}{#3}%
\cell{5.5}{2.5}{l}{#4}%
\cell{3}{4.2}{b}{#5}%
\cell{3}{0.8}{t}{#6}%
\cell{3.2}{2.5}{l}{#7}%
\qbezier(3.5,6.5)(7,2.5)(3.5,-1.5)%
\qbezier(2.5,6.5)(-1,2.5)(2.5,-1.5)%
\put(2.5,-1.5){\vector(1,-1){0}}%
\put(3.5,-1.5){\vector(-1,-1){0}}%
\qbezier(1,3)(3,5)(5,3)%
\qbezier(1,2)(3,0)(5,2)%
\put(5,3){\vector(1,-1){0}}%
\put(5,2){\vector(1,1){0}}%
\put(3,3.5){\vector(0,-1){2}}%
\end{picture}}

\mcm{\cthreecellv}{7}{%
\cinitdims{5}{8.2}%
\abovepic{#1}%
\belowpic{#2}%
\sidespic{#3}%
\sidespic{#4}%
\present{\precthreecellv{#1}{#2}{#3}{#4}{#5}{#6}{#7}}}

\newcommand{\pretopez}[2]{%
\begin{picture}(2.6,2.3)(-1.3,-2.2)%
\cell{0}{-2.2}{t}{#1}%
\cell{0}{-1.2}{c}{#2}%
\qbezier(0,0)(-2,-2)(0,-2)%
\qbezier(0,0)(2,-2)(0,-2)%
\put(0,0){\vector(-1,1){0}}%
\end{picture}}

\mcm{\topez}{2}{%
\ginitdims{2.6}{2.3}%
\belowpic{#1}%
\present{\pretopez{#1}{#2}}}

\newcommand{\pretopea}[3]{%
\begin{picture}(4,1.9)(-2,-0,2)%
\cell{0}{1.7}{b}{#1}%
\cell{0}{-0.2}{t}{#2}%
\cell{0}{0.7}{c}{#3}%
\qbezier(-2,0)(0,3)(2,0)%
\put(-2,0){\vector(1,0){4}}%
\put(2,0){\vector(2,-3){0}}%
\end{picture}}

\mcm{\topea}{3}{%
\ginitdims{4}{1.9}%
\abovepic{#1}%
\belowpic{#2}%
\present{\pretopea{#1}{#2}{#3}}}

\newcommand{\pretopeb}[4]{%
\begin{picture}(4,2.2)(-2,-0.2)%
\cell{-1.1}{1}{br}{#1}%
\cell{1.1}{1}{bl}{#2}%
\cell{0}{-0.2}{t}{#3}%
\cell{0}{0.8}{c}{#4}%
\put(-2,0){\vector(1,1){2}}%
\put(0,2){\vector(1,-1){2}}%
\put(-2,0){\vector(1,0){4}}%
\end{picture}}

\mcm{\topeb}{4}{%
\ginitdims{4}{2.2}%
\belowpic{#3}%
\present{\pretopeb{#1}{#2}{#3}{#4}}}

\newcommand{\pretopec}[5]{%
\begin{picture}(4,2.2)(-2,-0.2)%
\cell{-1.8}{1}{br}{#1}%
\cell{0}{2.2}{b}{#2}%
\cell{1.8}{1}{bl}{#3}%
\cell{0}{-0.2}{t}{#4}%
\cell{0}{0.8}{c}{#5}%
\put(-2,0){\vector(1,2){1}}%
\put(-1,2){\vector(1,0){2}}%
\put(1,2){\vector(1,-2){1}}%
\put(-2,0){\vector(1,0){4}}%
\end{picture}}

\mcm{\topec}{5}{%
\ginitdims{4}{2.2}%
\sidespic{#1}%
\abovepic{#2}%
\sidespic{#3}%
\belowpic{#4}%
\present{\pretopec{#1}{#2}{#3}{#4}{#5}}}

\newcommand{\pretoped}[6]{%
\begin{picture}(4,2.5)(-2,-0.2)%
\cell{-2}{0.6}{br}{#1}%
\cell{-0.7}{2.2}{br}{#2}%
\cell{0.7}{2.2}{bl}{#3}%
\cell{2}{0.6}{bl}{#4}%
\cell{0}{-0.2}{t}{#5}%
\cell{0}{0.8}{c}{#6}%
\put(-2,0){\vector(1,3){0.5}}%
\put(-1.5,1.5){\vector(3,2){1.5}}%
\put(0,2.5){\vector(3,-2){1.5}}%
\put(1.5,1.5){\vector(1,-3){0.5}}%
\put(-2,0){\vector(1,0){4}}%
\end{picture}}

\mcm{\toped}{6}{%
\ginitdims{4}{2.5}%
\sidespic{#1}%
\abovepic{#2}%
\abovepic{#3}%
\sidespic{#4}%
\belowpic{#5}%
\present{\pretoped{#1}{#2}{#3}{#4}{#5}{#6}}}

\newcommand{\pretopeq}[5]{%
\begin{picture}(4,2.5)(-2,-0.2)%
\cell{-2}{0.6}{br}{#1}%
\cell{-1}{2.2}{br}{#2}%
\cell{2}{0.6}{bl}{#3}%
\cell{0}{-0.2}{t}{#4}%
\cell{0}{0.8}{c}{#5}%
\put(-2,0){\vector(1,3){0.5}}%
\put(-1.5,1.5){\vector(1,1){1}}%
\cell{0.9}{2.3}{c}{\ddots}
\put(1.5,1.5){\vector(1,-3){0.5}}%
\put(-2,0){\vector(1,0){4}}%
\end{picture}}

\mcm{\topeq}{5}{%
\ginitdims{4}{2.5}%
\sidespic{#1}%
\abovepic{#2}%
\sidespic{#3}%
\belowpic{#4}%
\present{\pretopeq{#1}{#2}{#3}{#4}{#5}}}

\newcommand{\pretopebase}[1]{%
\begin{picture}(4,0.4)(0,-0.2)%
\cell{2}{0.2}{b}{#1}%
\put(0,0){\vector(1,0){4}}%
\end{picture}}

\mcm{\topebase}{1}{%
\ginitdims{4}{0.4}%
\abovepic{#1}%
\present{\pretopebase{#1}}}

\newcommand{\pretopezs}[2]{%
\begin{picture}(2.6,2.3)(-1.3,-2.2)%
\cell{0}{-2.2}{t}{#1}%
\cell{0}{-1.2}{c}{#2}%
\qbezier(0,0)(-2,-2)(0,-2)%
\qbezier(0,0)(2,-2)(0,-2)%
\end{picture}}

\mcm{\topezs}{2}{%
\ginitdims{2.6}{2.3}%
\belowpic{#1}%
\present{\pretopezs{#1}{#2}}}

\newcommand{\pretopeas}[3]{%
\begin{picture}(4,1.9)(-2,-0,2)%
\cell{0}{1.7}{b}{#1}%
\cell{0}{-0.2}{t}{#2}%
\cell{0}{0.7}{c}{#3}%
\qbezier(-2,0)(0,3)(2,0)%
\put(-2,0){\line(1,0){4}}%
\end{picture}}

\mcm{\topeas}{3}{%
\ginitdims{4}{1.9}%
\abovepic{#1}%
\belowpic{#2}%
\present{\pretopeas{#1}{#2}{#3}}}

\newcommand{\pretopebs}[4]{%
\begin{picture}(4,2.2)(-2,-0.2)%
\cell{-1.1}{1}{br}{#1}%
\cell{1.1}{1}{bl}{#2}%
\cell{0}{-0.2}{t}{#3}%
\cell{0}{0.8}{c}{#4}%
\put(-2,0){\line(1,1){2}}%
\put(0,2){\line(1,-1){2}}%
\put(-2,0){\line(1,0){4}}%
\end{picture}}

\mcm{\topebs}{4}{%
\ginitdims{4}{2.2}%
\belowpic{#3}%
\present{\pretopebs{#1}{#2}{#3}{#4}}}

\newcommand{\pretopecs}[5]{%
\begin{picture}(4,2.2)(-2,-0.2)%
\cell{-1.8}{1}{br}{#1}%
\cell{0}{2.2}{b}{#2}%
\cell{1.8}{1}{bl}{#3}%
\cell{0}{-0.2}{t}{#4}%
\cell{0}{0.8}{c}{#5}%
\put(-2,0){\line(1,2){1}}%
\put(-1,2){\line(1,0){2}}%
\put(1,2){\line(1,-2){1}}%
\put(-2,0){\line(1,0){4}}%
\end{picture}}

\mcm{\topecs}{5}{%
\ginitdims{4}{2.2}%
\sidespic{#1}%
\abovepic{#2}%
\sidespic{#3}%
\belowpic{#4}%
\present{\pretopecs{#1}{#2}{#3}{#4}{#5}}}

\newcommand{\pretopeds}[6]{%
\begin{picture}(4,2.5)(-2,-0.2)%
\cell{-2}{0.6}{br}{#1}%
\cell{-0.7}{2.2}{br}{#2}%
\cell{0.7}{2.2}{bl}{#3}%
\cell{2}{0.6}{bl}{#4}%
\cell{0}{-0.2}{t}{#5}%
\cell{0}{0.8}{c}{#6}%
\put(-2,0){\line(1,3){0.5}}%
\put(-1.5,1.5){\line(3,2){1.5}}%
\put(0,2.5){\line(3,-2){1.5}}%
\put(1.5,1.5){\line(1,-3){0.5}}%
\put(-2,0){\line(1,0){4}}%
\end{picture}}

\mcm{\topeds}{6}{%
\ginitdims{4}{2.5}%
\sidespic{#1}%
\abovepic{#2}%
\abovepic{#3}%
\sidespic{#4}%
\belowpic{#5}%
\present{\pretopeds{#1}{#2}{#3}{#4}{#5}{#6}}}

\newcommand{\pretopeqs}[5]{%
\begin{picture}(4,2.5)(-2,-0.2)%
\cell{-2}{0.6}{br}{#1}%
\cell{-1}{2.2}{br}{#2}%
\cell{2}{0.6}{bl}{#3}%
\cell{0}{-0.2}{t}{#4}%
\cell{0}{0.8}{c}{#5}%
\put(-2,0){\line(1,3){0.5}}%
\put(-1.5,1.5){\line(1,1){1}}%
\cell{0.9}{2.3}{c}{\ddots}
\put(1.5,1.5){\line(1,-3){0.5}}%
\put(-2,0){\line(1,0){4}}%
\end{picture}}

\mcm{\topeqs}{5}{%
\ginitdims{4}{2.5}%
\sidespic{#1}%
\abovepic{#2}%
\sidespic{#3}%
\belowpic{#4}%
\present{\pretopeqs{#1}{#2}{#3}{#4}{#5}}}

\newcommand{\pretopebases}[1]{%
\begin{picture}(4,0.4)(0,-0.2)%
\cell{2}{0.2}{b}{#1}%
\put(0,0){\line(1,0){4}}%
\end{picture}}

\mcm{\topebases}{1}{%
\ginitdims{4}{0.4}%
\abovepic{#1}%
\present{\pretopebases{#1}}}

\newcommand{\pregdots}[6]{%
\begin{picture}(5,8.4)(0,-2.7)%
\cell{2.5}{5.7}{b}{#1}%
\cell{1.5}{2.8}{b}{#2}%
\cell{1.5}{0.2}{t}{#3}%
\cell{2.5}{-2.7}{t}{#4}%
\cell{2.7}{4.25}{l}{#5}%
\cell{2.7}{-1.25}{l}{#6}%
\qbezier(0,1.5)(2.5,9.5)(5,1.5)%
\qbezier(0,1.5)(2.5,4)(5,1.5)%
\qbezier(0,1.5)(2.5,-1)(5,1.5)%
\qbezier(0,1.5)(2.5,-6.5)(5,1.5)%
\put(2.5,5.25){\vector(0,-1){2}}%
\put(2.5,-0.25){\vector(0,-1){2}}%
\cell{2.5}{1.7}{c}{\vdots}%
\put(5,1.5){\vector(1,-4){0}}%
\put(5,1.5){\vector(4,-3){0}}%
\put(5,1.5){\vector(4,3){0}}%
\put(5,1.5){\vector(1,4){0}}%
\end{picture}}

\mcm{\gdots}{6}{%
\ginitdims{5}{8.4}%
\abovepic{#1}%
\belowpic{#4}%
\present{\pregdots{#1}{#2}{#3}{#4}{#5}{#6}}}

\newlength{\volt}
\makeatletter

\def\diagram{\m@th\leftwidth=\z@ \rightwidth=\z@ \topheight=\z@
\botheight=\z@ \setbox\@picbox\hbox\bgroup}

\def\enddiagram{\egroup\wd\@picbox\rightwidth\unitlength
\ht\@picbox\topheight\unitlength \dp\@picbox\botheight\unitlength
\hskip\leftwidth\unitlength\box\@picbox}

\def\bfig{\begin{diagram}}
\def\efig{\end{diagram}}
\newcount\wideness \newcount\leftwidth \newcount\rightwidth
\newcount\highness \newcount\topheight \newcount\botheight

\def\ratchet#1#2{\ifnum#1<#2 \global #1=#2 \fi}

\def\putbox(#1,#2)#3{%
\horsize{\wideness}{#3} \divide\wideness by 2 {\advance\wideness
by #1 \ratchet{\rightwidth}{\wideness}} {\advance\wideness by -#1
\ratchet{\leftwidth}{\wideness}} \vertsize{\highness}{#3}
\divide\highness by 2 {\advance\highness by #2
\ratchet{\topheight}{\highness}} {\advance\highness by -#2
\ratchet{\botheight}{\highness}} \put(#1,#2){\makebox(0,0){$#3$}}}

\def\putlbox(#1,#2)#3{%
\horsize{\wideness}{#3} {\advance\wideness by #1
\ratchet{\rightwidth}{\wideness}} {\ratchet{\leftwidth}{-#1}}
\vertsize{\highness}{#3} \divide\highness by 2 {\advance\highness
by #2 \ratchet{\topheight}{\highness}} {\advance\highness by -#2
\ratchet{\botheight}{\highness}}
\put(#1,#2){\makebox(0,0)[l]{$#3$}}}

\def\putrbox(#1,#2)#3{%
\horsize{\wideness}{#3} {\ratchet{\rightwidth}{#1}}
{\advance\wideness by -#1 \ratchet{\leftwidth}{\wideness}}
\vertsize{\highness}{#3} \divide\highness by 2 {\advance\highness
by #2 \ratchet{\topheight}{\highness}} {\advance\highness by -#2
\ratchet{\botheight}{\highness}}
\put(#1,#2){\makebox(0,0)[r]{$#3$}}}

\def\adjust[#1]{} 

\newcount \coefa
\newcount \coefb
\newcount \coefc
\newcount\tempcounta
\newcount\tempcountb
\newcount\tempcountc
\newcount\tempcountd
\newcount\xext
\newcount\yext
\newcount\xoff
\newcount\yoff
\newcount\gap%
\newcount\arrowtypea
\newcount\arrowtypeb
\newcount\arrowtypec
\newcount\arrowtyped
\newcount\arrowtypee
\newcount\height
\newcount\width
\newcount\xpos
\newcount\ypos
\newcount\run
\newcount\rise
\newcount\arrowlength
\newcount\halflength
\newcount\arrowtype
\newdimen\tempdimen
\newdimen\xlen
\newdimen\ylen
\newsavebox{\tempboxa}%
\newsavebox{\tempboxb}%
\newsavebox{\tempboxc}%

\newdimen\w@dth

\def\setw@dth#1#2{\setbox\z@\hbox{\m@th$#1$}\w@dth=\wd\z@
\setbox\@ne\hbox{\m@th$#2$}\ifnum\w@dth<\wd\@ne \w@dth=\wd\@ne \fi
\advance\w@dth by 1.2em}

\def\t@^#1_#2{\allowbreak\def\n@one{#1}\def\n@two{#2}\mathrel
{\setw@dth{#1}{#2} \mathop{\hbox to
\w@dth{\rightarrowfill}}\limits \ifx\n@one\empty\else
^{\box\z@}\fi \ifx\n@two\empty\else _{\box\@ne}\fi}}
\def\t@@^#1{\@ifnextchar_{\t@^{#1}}{\t@^{#1}_{}}}
\def\to{\@ifnextchar^{\t@@}{\t@@^{}}}

\def\t@left^#1_#2{\def\n@one{#1}\def\n@two{#2}\mathrel{\setw@dth{#1}{#2}
\mathop{\hbox to \w@dth{\leftarrowfill}}\limits
\ifx\n@one\empty\else ^{\box\z@}\fi \ifx\n@two\empty\else
_{\box\@ne}\fi}}
\def\t@@left^#1{\@ifnextchar_{\t@left^{#1}}{\t@left^{#1}_{}}}
\def\toleft{\@ifnextchar^{\t@@left}{\t@@left^{}}}

\def\two@^#1_#2{\allowbreak
\def\n@one{#1}\def\n@two{#2}\mathrel{\setw@dth{#1}{#2}
\mathop{\vcenter{\lineskip\z@\baselineskip\z@
                 \hbox to \w@dth{\rightarrowfill}%
                 \hbox to \w@dth{\rightarrowfill}}%
       }\limits
\ifx\n@one\empty\else ^{\box\z@}\fi \ifx\n@two\empty\else
_{\box\@ne}\fi}}
\def\tw@@^#1{\@ifnextchar _{\two@^{#1}}{\two@^{#1}_{}}}
\def\two{\@ifnextchar ^{\tw@@}{\tw@@^{}}}

\def\tofr@^#1_#2{\def\n@one{#1}\def\n@two{#2}\mathrel{\setw@dth{#1}{#2}
\mathop{\vcenter{\hbox to \w@dth{\rightarrowfill}\kern-1.7ex
                 \hbox to \w@dth{\leftarrowfill}}%
       }\limits
\ifx\n@one\empty\else ^{\box\z@}\fi \ifx\n@two\empty\else
_{\box\@ne}\fi}}
\def\t@fr@^#1{\@ifnextchar_ {\tofr@^{#1}}{\tofr@^{#1}_{}}}
\def\tofro{\@ifnextchar^ {\t@fr@}{\t@fr@^{}}}

\def\mon{\mathop{\m@th\hbox to
      14.6\P@{\lasyb\char'51\hskip-2.1\P@$\arrext$\hss
$\mathord\rightarrow$}}\limits} 
\def\leftmono{\mathrel{\m@th\hbox to
14.6\P@{$\mathord\leftarrow$\hss$\arrext$\hskip-2.1\P@\lasyb\char'50%
}}\limits} 
\mathchardef\arrext="0200       

\def\settypes(#1,#2,#3){\arrowtypea#1 \arrowtypeb#2 \arrowtypec#3}
\def\settoheight#1#2{\setbox\@tempboxa\hbox{#2}#1\ht\@tempboxa\relax}%
\def\settodepth#1#2{\setbox\@tempboxa\hbox{#2}#1\dp\@tempboxa\relax}%
\def\settokens`#1`#2`#3`#4`{%
     \def\tokena{#1}\def\tokenb{#2}\def\tokenc{#3}\def\tokend{#4}}
\def\setsqparms[#1`#2`#3`#4;#5`#6]{%
\arrowtypea #1 \arrowtypeb #2 \arrowtypec #3 \arrowtyped #4
\width #5 \height #6 }
\def\setpos(#1,#2){\xpos=#1 \ypos#2}

\def\settriparms[#1`#2`#3;#4]{\settripairparms[#1`#2`#3`1`1;#4]}%

\def\settripairparms[#1`#2`#3`#4`#5;#6]{%
\arrowtypea #1 \arrowtypeb #2 \arrowtypec #3 \arrowtyped #4
\arrowtypee #5 \width #6 \height #6 }

\def\resetparms{\settripairparms[1`1`1`1`1;500]\width 500}

\resetparms

\def\mvector(#1,#2)#3{
\put(0,0){\vector(#1,#2){#3}}%
\put(0,0){\vector(#1,#2){26}}%
}
\def\evector(#1,#2)#3{{
\arrowlength #3
\put(0,0){\vector(#1,#2){\arrowlength}}%
\advance \arrowlength by-30
\put(0,0){\vector(#1,#2){\arrowlength}}%
}}

\def\horsize#1#2{%
\settowidth{\tempdimen}{$#2$}%
#1=\tempdimen \divide #1 by\unitlength }

\def\vertsize#1#2{%
\settoheight{\tempdimen}{$#2$}%
#1=\tempdimen
\settodepth{\tempdimen}{$#2$}%
\advance #1 by\tempdimen \divide #1 by\unitlength }

\def\putvector(#1,#2)(#3,#4)#5#6{{%
\ifnum3<\arrowtype \putdashvector(#1,#2)(#3,#4)#5\arrowtype \else
\ifnum\arrowtype<-3 \putdashvector(#1,#2)(#3,#4)#5\arrowtype \else
\xpos=#1 \ypos=#2 \run=#3 \rise=#4 \arrowlength=#5 \ifnum
\arrowtype<0
    \ifnum \run=0
        \advance \ypos by-\arrowlength
    \else
        \tempcounta \arrowlength
        \multiply \tempcounta by\rise
        \divide \tempcounta by\run
        \ifnum\run>0
            \advance \xpos by\arrowlength
            \advance \ypos by\tempcounta
        \else
            \advance \xpos by-\arrowlength
            \advance \ypos by-\tempcounta
        \fi
    \fi
    \multiply \arrowtype by-1
    \multiply \rise by-1
    \multiply \run by-1
\fi \ifcase \arrowtype
\or \put(\xpos,\ypos){\vector(\run,\rise){\arrowlength}}%
\or \put(\xpos,\ypos){\mvector(\run,\rise)\arrowlength}%
\or \put(\xpos,\ypos){\evector(\run,\rise){\arrowlength}}%
\fi\fi\fi }}

\def\putsplitvector(#1,#2)#3#4{
\xpos #1 \ypos #2 \arrowtype #4 \halflength #3 \arrowlength #3
\gap 140 \advance \halflength by-\gap \divide \halflength by2
\ifnum\arrowtype>0
   \ifcase \arrowtype
   \or \put(\xpos,\ypos){\line(0,-1){\halflength}}%
       \advance\ypos by-\halflength
       \advance\ypos by-\gap
       \put(\xpos,\ypos){\vector(0,-1){\halflength}}%
   \or \put(\xpos,\ypos){\line(0,-1)\halflength}%
       \put(\xpos,\ypos){\vector(0,-1)3}%
       \advance\ypos by-\halflength
       \advance\ypos by-\gap
       \put(\xpos,\ypos){\vector(0,-1){\halflength}}%
   \or \put(\xpos,\ypos){\line(0,-1)\halflength}%
       \advance\ypos by-\halflength
       \advance\ypos by-\gap
       \put(\xpos,\ypos){\evector(0,-1){\halflength}}%
   \fi
\else \arrowtype=-\arrowtype
   \ifcase\arrowtype
   \or \advance \ypos by-\arrowlength
       \put(\xpos,\ypos){\line(0,1){\halflength}}%
       \advance\ypos by\halflength
       \advance\ypos by\gap
       \put(\xpos,\ypos){\vector(0,1){\halflength}}%
   \or \advance \ypos by-\arrowlength
       \put(\xpos,\ypos){\line(0,1)\halflength}%
       \put(\xpos,\ypos){\vector(0,1)3}%
       \advance\ypos by\halflength
       \advance\ypos by\gap
       \put(\xpos,\ypos){\vector(0,1){\halflength}}%
   \or \advance \ypos by-\arrowlength
       \put(\xpos,\ypos){\line(0,1)\halflength}%
       \advance\ypos by\halflength
       \advance\ypos by\gap
       \put(\xpos,\ypos){\evector(0,1){\halflength}}%
   \fi
\fi }

\def\putmorphism(#1)(#2,#3)[#4`#5`#6]#7#8#9{{%
\run #2 \rise #3 \ifnum\rise=0
  \puthmorphism(#1)[#4`#5`#6]{#7}{#8}#9%
\else\ifnum\run=0
  \putvmorphism(#1)[#4`#5`#6]{#7}{#8}#9%
\else
\setpos(#1)%
\arrowlength #7 \arrowtype #8 \ifnum\run=0 \else\ifnum\rise=0
\else \ifnum\run>0
    \coefa=1
\else
   \coefa=-1
\fi \ifnum\arrowtype>0
   \coefb=0
   \coefc=-1
\else
   \coefb=\coefa
   \coefc=1
   \arrowtype=-\arrowtype
\fi \width=2 \multiply \width by\run \divide \width by\rise
\ifnum \width<0  \width=-\width\fi \advance\width by60 \if l#9
\width=-\width\fi
\putbox(\xpos,\ypos){#4}
{\multiply \coefa by\arrowlength
\advance\xpos by\coefa \multiply \coefa by\rise \divide \coefa
by\run \advance \ypos by\coefa
\putbox(\xpos,\ypos){#5} }%
{\multiply \coefa by\arrowlength
\divide \coefa by2 \advance \xpos by\coefa \advance \xpos by\width
\multiply \coefa by\rise \divide \coefa by\run \advance \ypos
by\coefa
\if l#9%
   \putrbox(\xpos,\ypos){#6}%
\else\if r#9%
   \putlbox(\xpos,\ypos){#6}%
\fi\fi }%
{\multiply \rise by-\coefc
\multiply \run by-\coefc \multiply \coefb by\arrowlength \advance
\xpos by\coefb \multiply \coefb by\rise \divide \coefb by\run
\advance \ypos by\coefb \multiply \coefc by70 \advance \ypos
by\coefc \multiply \coefc by\run \divide \coefc by\rise \advance
\xpos by\coefc \multiply \coefa by140 \multiply \coefa by\run
\divide \coefa by\rise \advance \arrowlength by\coefa
\ifcase\arrowtype
\or \put(\xpos,\ypos){\vector(\run,\rise){\arrowlength}}%
\or \put(\xpos,\ypos){\mvector(\run,\rise){\arrowlength}}%
\or \put(\xpos,\ypos){\evector(\run,\rise){\arrowlength}}%
\fi}\fi\fi\fi\fi}}

\newcount\numbdashes \newcount\lengthdash \newcount\increment

\def\howmanydashes{
\numbdashes=\arrowlength \lengthdash=40 \divide\numbdashes by
\lengthdash \lengthdash=\arrowlength \divide\lengthdash by
\numbdashes
\increment=\lengthdash \multiply\lengthdash by 3
\divide\lengthdash by 5 }

\def\putdashvector(#1)(#2,#3)#4#5{%
\ifnum#3=0 \putdashhvector(#1){#4}#5 \else \ifnum#2=0
\putdashvvector(#1){#4}#5\fi\fi}

\def\putdashhvector(#1,#2)#3#4{{%
\arrowlength=#3 \howmanydashes
\multiput(#1,#2)(\increment,0){\numbdashes}%
{\vrule height .4pt width \lengthdash\unitlength} \arrowtype=#4
\xpos=#1 \ifnum\arrowtype<0 \advance\arrowtype by 7 \fi
\ifcase\arrowtype \or \advance\xpos by 10
    \put(\xpos,#2){\vector(-1,0){\lengthdash}}
    \advance\xpos by 40
    \put(\xpos,#2){\vector(-1,0){\lengthdash}}
\or \advance \xpos by 10
    \put(\xpos,#2){\vector(-1,0){\lengthdash}}
    \advance\xpos by  \arrowlength
    \advance\xpos by  -50
    \put(\xpos,#2){\vector(-1,0){\lengthdash}}
\or \advance\xpos by 10
    \put(\xpos,#2){\vector(-1,0){\lengthdash}}
\or \advance\xpos by \arrowlength
    \advance\xpos by -\lengthdash
    \put(\xpos,#2){\vector(1,0){\lengthdash}}
\or {\advance\xpos by 10
    \put(\xpos,#2){\vector(1,0){\lengthdash}}}
    \advance\xpos by \arrowlength
    \advance\xpos by -\lengthdash
    \put(\xpos,#2){\vector(1,0){\lengthdash}}
\or \advance\xpos by \arrowlength
    \advance\xpos by -\lengthdash
    \put(\xpos,#2){\vector(1,0){\lengthdash}}
    \advance\xpos by -40
    \put(\xpos,#2){\vector(1,0){\lengthdash}}
   \fi
}}

\def\putdashvvector(#1,#2)#3#4{{%
\arrowlength=#3 \howmanydashes \ypos=#2 \advance\ypos by
-\arrowlength
\multiput(#1,#2)(0,\increment){\numbdashes}%
    {\vrule width .4pt height \lengthdash\unitlength}
\arrowtype=#4 \ypos=#2 \ifnum\arrowtype<0 \advance\arrowtype by 7
\fi \ifcase\arrowtype \or \advance\ypos by \arrowlength
\advance\ypos by -40
    \put(#1,\ypos){\vector(0,1){\lengthdash}}
    \advance\ypos by -40
    \put(#1,\ypos){\vector(0,1){\lengthdash}}
\or \advance\ypos by 10
    \put(#1,\ypos){\vector(0,1){\lengthdash}}
    \advance\ypos by \arrowlength \advance\ypos by -40
    \put(#1,\ypos){\vector(0,1){\lengthdash}}
\or \advance\ypos by \arrowlength \advance\ypos by -40
    \put(#1,\ypos){\vector(0,1){\lengthdash}}
\or \advance\ypos by 10
    \put(#1,\ypos){\vector(0,-1){\lengthdash}}
\or \advance\ypos by 10
    \put(#1,\ypos){\vector(0,-1){\lengthdash}}
    \advance\ypos by \arrowlength \advance\ypos by -40
    \put(#1,\ypos){\vector(0,-1){\lengthdash}}
\or \advance\ypos by 10
    \put(#1,\ypos){\vector(0,-1){\lengthdash}}
    \advance\ypos by 40
    \put(#1,\ypos){\vector(0,-1){\lengthdash}}
\fi }}

\def\puthmorphism(#1,#2)[#3`#4`#5]#6#7#8{{%
\xpos #1 \ypos #2 \width #6 \arrowlength #6 \arrowtype=#7
\putbox(\xpos,\ypos){#3\vphantom{#4}}%
{\advance \xpos by\arrowlength
\putbox(\xpos,\ypos){\vphantom{#3}#4}}%
\horsize{\tempcounta}{#3}%
\horsize{\tempcountb}{#4}%
\divide \tempcounta by2 \divide \tempcountb by2 \advance
\tempcounta by30 \advance \tempcountb by30 \advance \xpos
by\tempcounta \advance \arrowlength by-\tempcounta \advance
\arrowlength by-\tempcountb
\putvector(\xpos,\ypos)(1,0)\arrowlength\arrowtype \divide
\arrowlength by2 \advance \xpos by\arrowlength
\vertsize{\tempcounta}{#5}%
\divide\tempcounta by2 \advance \tempcounta by20
\if a#8 %
   \advance \ypos by\tempcounta
   \putbox(\xpos,\ypos){#5}%
\else
   \advance \ypos by-\tempcounta
   \putbox(\xpos,\ypos){#5}%
\fi}}

\def\putvmorphism(#1,#2)[#3`#4`#5]#6#7#8{{%
\xpos #1 \ypos #2 \arrowlength #6 \arrowtype #7
\settowidth{\xlen}{$#5$}%
\putbox(\xpos,\ypos){#3}%
{\advance \ypos by-\arrowlength
\putbox(\xpos,\ypos){#4}}%
{\advance\arrowlength by-140 \advance \ypos by-70 \ifdim\xlen>0pt
   \if m#8%
      \putsplitvector(\xpos,\ypos)\arrowlength\arrowtype
   \else
   \putvector(\xpos,\ypos)(0,-1)\arrowlength\arrowtype
   \fi
\else
   \putvector(\xpos,\ypos)(0,-1)\arrowlength\arrowtype
\fi}%
\ifdim\xlen>0pt
   \divide \arrowlength by2
   \advance\ypos by-\arrowlength
   \if l#8%
      \advance \xpos by-40
      \putrbox(\xpos,\ypos){#5}%
   \else\if r#8%
      \advance \xpos by40
      \putlbox(\xpos,\ypos){#5}%
   \else
      \putbox(\xpos,\ypos){#5}%
   \fi\fi
\fi }}

\def\putsquarep<#1>(#2)[#3;#4`#5`#6`#7]{{%
\setsqparms[#1]%
\setpos(#2)%
\settokens`#3`%
\puthmorphism(\xpos,\ypos)[\tokenc`\tokend`{#7}]{\width}{\arrowtyped}b%
\advance\ypos by \height
\puthmorphism(\xpos,\ypos)[\tokena`\tokenb`{#4}]{\width}{\arrowtypea}a%
\putvmorphism(\xpos,\ypos)[``{#5}]{\height}{\arrowtypeb}l%
\advance\xpos by \width
\putvmorphism(\xpos,\ypos)[``{#6}]{\height}{\arrowtypec}r%
}}

\def\putsquare{\@ifnextchar <{\putsquarep}{\putsquarep%
   <\arrowtypea`\arrowtypeb`\arrowtypec`\arrowtyped;\width`\height>}}
\def\square{\@ifnextchar< {\squarep}{\squarep
   <\arrowtypea`\arrowtypeb`\arrowtypec`\arrowtyped;\width`\height>}}
\def\squarep<#1>[#2`#3`#4`#5;#6`#7`#8`#9]{{
\setsqparms[#1]
\diagram
\putsquarep<\arrowtypea`\arrowtypeb`\arrowtypec`
\arrowtyped;\width`\height>
(0,0)[#2`#3`#4`{#5};#6`#7`#8`{#9}]
\enddiagram
}}                                                 
\def\putptrianglep<#1>(#2,#3)[#4`#5`#6;#7`#8`#9]{{%
\settriparms[#1]%
\xpos=#2 \ypos=#3 \advance\ypos by \height
\puthmorphism(\xpos,\ypos)[#4`#5`{#7}]{\height}{\arrowtypea}a%
\putvmorphism(\xpos,\ypos)[`#6`{#8}]{\height}{\arrowtypeb}l%
\advance\xpos by\height
\putmorphism(\xpos,\ypos)(-1,-1)[``{#9}]{\height}{\arrowtypec}r%
}}

\def\putptriangle{\@ifnextchar <{\putptrianglep}{\putptrianglep
   <\arrowtypea`\arrowtypeb`\arrowtypec;\height>}}
\def\ptriangle{\@ifnextchar <{\ptrianglep}{\ptrianglep
   <\arrowtypea`\arrowtypeb`\arrowtypec;\height>}}
\def\ptrianglep<#1>[#2`#3`#4;#5`#6`#7]{{
\settriparms[#1]
\diagram
\putptrianglep<\arrowtypea`\arrowtypeb`
\arrowtypec;\height>
(0,0)[#2`#3`#4;#5`#6`{#7}]
\enddiagram
}}                                            

\def\putqtrianglep<#1>(#2,#3)[#4`#5`#6;#7`#8`#9]{{%
\settriparms[#1]%
\xpos=#2 \ypos=#3 \advance\ypos by\height
\puthmorphism(\xpos,\ypos)[#4`#5`{#7}]{\height}{\arrowtypea}a%
\putmorphism(\xpos,\ypos)(1,-1)[``{#8}]{\height}{\arrowtypeb}l%
\advance\xpos by\height
\putvmorphism(\xpos,\ypos)[`#6`{#9}]{\height}{\arrowtypec}r%
}}

\def\putqtriangle{\@ifnextchar <{\putqtrianglep}{\putqtrianglep
   <\arrowtypea`\arrowtypeb`\arrowtypec;\height>}}
\def\qtriangle{\@ifnextchar <{\qtrianglep}{\qtrianglep
   <\arrowtypea`\arrowtypeb`\arrowtypec;\height>}}
\def\qtrianglep<#1>[#2`#3`#4;#5`#6`#7]{{
\settriparms[#1]
\width=\height                                
\diagram
\putqtrianglep<\arrowtypea`\arrowtypeb`
\arrowtypec;\height>
(0,0)[#2`#3`#4;#5`#6`{#7}]
\enddiagram
}}

\def\putdtrianglep<#1>(#2,#3)[#4`#5`#6;#7`#8`#9]{{%
\settriparms[#1]%
\xpos=#2 \ypos=#3
\puthmorphism(\xpos,\ypos)[#5`#6`{#9}]{\height}{\arrowtypec}b%
\advance\xpos by \height \advance\ypos by\height
\putmorphism(\xpos,\ypos)(-1,-1)[``{#7}]{\height}{\arrowtypea}l%
\putvmorphism(\xpos,\ypos)[#4``{#8}]{\height}{\arrowtypeb}r%
}}

\def\putdtriangle{\@ifnextchar <{\putdtrianglep}{\putdtrianglep
   <\arrowtypea`\arrowtypeb`\arrowtypec;\height>}}
\def\dtriangle{\@ifnextchar <{\dtrianglep}{\dtrianglep
   <\arrowtypea`\arrowtypeb`\arrowtypec;\height>}}
\def\dtrianglep<#1>[#2`#3`#4;#5`#6`#7]{{
\settriparms[#1]
\width=\height                                
\diagram
\putdtrianglep<\arrowtypea`\arrowtypeb`
\arrowtypec;\height>
(0,0)[#2`#3`#4;#5`#6`{#7}]
\enddiagram
}}

\def\putbtrianglep<#1>(#2,#3)[#4`#5`#6;#7`#8`#9]{{%
\settriparms[#1]%
\xpos=#2 \ypos=#3
\puthmorphism(\xpos,\ypos)[#5`#6`{#9}]{\height}{\arrowtypec}b%
\advance\ypos by\height
\putmorphism(\xpos,\ypos)(1,-1)[``{#8}]{\height}{\arrowtypeb}r%
\putvmorphism(\xpos,\ypos)[#4``{#7}]{\height}{\arrowtypea}l%
}}

\def\putbtriangle{\@ifnextchar <{\putbtrianglep}{\putbtrianglep
   <\arrowtypea`\arrowtypeb`\arrowtypec;\height>}}
\def\btriangle{\@ifnextchar <{\btrianglep}{\btrianglep
   <\arrowtypea`\arrowtypeb`\arrowtypec;\height>}}
\def\btrianglep<#1>[#2`#3`#4;#5`#6`#7]{{
\settriparms[#1]
\width=\height                               
\diagram
\putbtrianglep<\arrowtypea`\arrowtypeb`
\arrowtypec;\height>
(0,0)[#2`#3`#4;#5`#6`{#7}]
\enddiagram
}}

\def\putAtrianglep<#1>(#2,#3)[#4`#5`#6;#7`#8`#9]{{%
\settriparms[#1]%
\xpos=#2 \ypos=#3 {\multiply \height by2
\puthmorphism(\xpos,\ypos)[#5`#6`{#9}]{\height}{\arrowtypec}b}%
\advance\xpos by\height \advance\ypos by\height
\putmorphism(\xpos,\ypos)(-1,-1)[#4``{#7}]{\height}{\arrowtypea}l%
\putmorphism(\xpos,\ypos)(1,-1)[``{#8}]{\height}{\arrowtypeb}r%
}}

\def\putAtriangle{\@ifnextchar <{\putAtrianglep}{\putAtrianglep
   <\arrowtypea`\arrowtypeb`\arrowtypec;\height>}}
\def\Atriangle{\@ifnextchar <{\Atrianglep}{\Atrianglep
   <\arrowtypea`\arrowtypeb`\arrowtypec;\height>}}
\def\Atrianglep<#1>[#2`#3`#4;#5`#6`#7]{{
\settriparms[#1]
\width=\height                                     
\diagram
\putAtrianglep<\arrowtypea`\arrowtypeb`
\arrowtypec;\height>
(0,0)[#2`#3`#4;#5`#6`{#7}]
\enddiagram
}}

\def\putAtrianglepairp<#1>(#2)[#3;#4`#5`#6`#7`#8]{{%
\settripairparms[#1]%
\setpos(#2)%
\settokens`#3`%
\puthmorphism(\xpos,\ypos)[\tokenb`\tokenc`{#7}]{\height}{\arrowtyped}b%
\advance\xpos by\height
\puthmorphism(\xpos,\ypos)[\phantom{\tokenc}`\tokend`{#8}]%
{\height}{\arrowtypee}b%
\advance\ypos by\height
\putmorphism(\xpos,\ypos)(-1,-1)[\tokena``{#4}]{\height}{\arrowtypea}l%
\putvmorphism(\xpos,\ypos)[``{#5}]{\height}{\arrowtypeb}m%
\putmorphism(\xpos,\ypos)(1,-1)[``{#6}]{\height}{\arrowtypec}r%
}}

\def\putAtrianglepair{\@ifnextchar <{\putAtrianglepairp}{\putAtrianglepairp%
   <\arrowtypea`\arrowtypeb`\arrowtypec`\arrowtyped`\arrowtypee;\height>}}
\def\Atrianglepair{\@ifnextchar <{\Atrianglepairp}{\Atrianglepairp%
   <\arrowtypea`\arrowtypeb`\arrowtypec`\arrowtyped`\arrowtypee;\height>}}

\def\Atrianglepairp<#1>[#2;#3`#4`#5`#6`#7]{{
\settripairparms[#1]
\settokens`#2`
\width=\height                                
\diagram
\putAtrianglepairp                            
<\arrowtypea`\arrowtypeb`\arrowtypec`
\arrowtyped`\arrowtypee;\height>
(0,0)[{#2};#3`#4`#5`#6`{#7}]
\enddiagram
}}

\def\putVtrianglep<#1>(#2,#3)[#4`#5`#6;#7`#8`#9]{{%
\settriparms[#1]%
\xpos=#2 \ypos=#3 \advance\ypos by\height {\multiply\height by2
\puthmorphism(\xpos,\ypos)[#4`#5`{#7}]{\height}{\arrowtypea}a}%
\putmorphism(\xpos,\ypos)(1,-1)[`#6`{#8}]{\height}{\arrowtypeb}l%
\advance\xpos by\height \advance\xpos by\height
\putmorphism(\xpos,\ypos)(-1,-1)[``{#9}]{\height}{\arrowtypec}r%
}}

\def\putVtriangle{\@ifnextchar <{\putVtrianglep}{\putVtrianglep
   <\arrowtypea`\arrowtypeb`\arrowtypec;\height>}}
\def\Vtriangle{\@ifnextchar <{\Vtrianglep}{\Vtrianglep
   <\arrowtypea`\arrowtypeb`\arrowtypec;\height>}}
\def\Vtrianglep<#1>[#2`#3`#4;#5`#6`#7]{{
\settriparms[#1]
\width=\height                                 
\diagram
\putVtrianglep<\arrowtypea`\arrowtypeb`
\arrowtypec;\height>
(0,0)[#2`#3`#4;#5`#6`{#7}]
\enddiagram
}}

\def\putVtrianglepairp<#1>(#2)[#3;#4`#5`#6`#7`#8]{{
\settripairparms[#1]%
\setpos(#2)%
\settokens`#3`%
\advance\ypos by\height
\putmorphism(\xpos,\ypos)(1,-1)[`\tokend`{#6}]{\height}{\arrowtypec}l%
\puthmorphism(\xpos,\ypos)[\tokena`\tokenb`{#4}]{\height}{\arrowtypea}a%
\advance\xpos by\height
\puthmorphism(\xpos,\ypos)[\phantom{\tokenb}`\tokenc`{#5}]%
{\height}{\arrowtypeb}a%
\putvmorphism(\xpos,\ypos)[``{#7}]{\height}{\arrowtyped}m%
\advance\xpos by\height
\putmorphism(\xpos,\ypos)(-1,-1)[``{#8}]{\height}{\arrowtypee}r%
}}

\def\putVtrianglepair{\@ifnextchar <{\putVtrianglepairp}{\putVtrianglepairp%
    <\arrowtypea`\arrowtypeb`\arrowtypec`\arrowtyped`\arrowtypee;\height>}}
\def\Vtrianglepair{\@ifnextchar <{\Vtrianglepairp}{\Vtrianglepairp%
    <\arrowtypea`\arrowtypeb`\arrowtypec`\arrowtyped`\arrowtypee;\height>}}
\def\Vtrianglepairp<#1>[#2;#3`#4`#5`#6`#7]{{
\settripairparms[#1]
\settokens`#2`
\diagram
\putVtrianglepairp                             
<\arrowtypea`\arrowtypeb`\arrowtypec`
\arrowtyped`\arrowtypee;\height>
(0,0)[{#2};#3`#4`#5`#6`{#7}]
\enddiagram
}}

\def\putCtrianglep<#1>(#2,#3)[#4`#5`#6;#7`#8`#9]{{%
\settriparms[#1]%
\xpos=#2 \ypos=#3 \advance\ypos by\height
\putmorphism(\xpos,\ypos)(1,-1)[``{#9}]{\height}{\arrowtypec}l%
\advance\xpos by\height \advance\ypos by\height
\putmorphism(\xpos,\ypos)(-1,-1)[#4`#5`{#7}]{\height}{\arrowtypea}l%
{\multiply\height by 2
\putvmorphism(\xpos,\ypos)[`#6`{#8}]{\height}{\arrowtypeb}r}%
}}

\def\putCtriangle{\@ifnextchar <{\putCtrianglep}{\putCtrianglep
    <\arrowtypea`\arrowtypeb`\arrowtypec;\height>}}
\def\Ctriangle{\@ifnextchar <{\Ctrianglep}{\Ctrianglep
    <\arrowtypea`\arrowtypeb`\arrowtypec;\height>}}
\def\Ctrianglep<#1>[#2`#3`#4;#5`#6`#7]{{
\settriparms[#1]
\width=\height                               
\diagram
\putCtrianglep<\arrowtypea`\arrowtypeb`
\arrowtypec;\height>
(0,0)[#2`#3`#4;#5`#6`{#7}]
\enddiagram
}}                                           
\def\putDtrianglep<#1>(#2,#3)[#4`#5`#6;#7`#8`#9]{{%
\settriparms[#1]%
\xpos=#2 \ypos=#3 \advance\xpos by\height \advance\ypos by\height
\putmorphism(\xpos,\ypos)(-1,-1)[``{#9}]{\height}{\arrowtypec}r%
\advance\xpos by-\height \advance\ypos by\height
\putmorphism(\xpos,\ypos)(1,-1)[`#5`{#8}]{\height}{\arrowtypeb}r%
{\multiply\height by 2
\putvmorphism(\xpos,\ypos)[#4`#6`{#7}]{\height}{\arrowtypea}l}%
}}

\def\putDtriangle{\@ifnextchar <{\putDtrianglep}{\putDtrianglep
    <\arrowtypea`\arrowtypeb`\arrowtypec;\height>}}
\def\Dtriangle{\@ifnextchar <{\Dtrianglep}{\Dtrianglep
   <\arrowtypea`\arrowtypeb`\arrowtypec;\height>}}
\def\Dtrianglep<#1>[#2`#3`#4;#5`#6`#7]{{
\settriparms[#1]
\width=\height                              
\diagram
\putDtrianglep<\arrowtypea`\arrowtypeb`
\arrowtypec;\height>
(0,0)[#2`#3`#4;#5`#6`{#7}]
\enddiagram
}}                                          
\def\setrecparms[#1`#2]{\width=#1 \height=#2}%

\def\recursep<#1`#2>[#3;#4`#5`#6`#7`#8]{{\m@th
\width=#1 \height=#2 \settokens`#3`
\settowidth{\tempdimen}{$\tokena$} \ifdim\tempdimen=0pt
  \savebox{\tempboxa}{\hbox{$\tokenb$}}%
  \savebox{\tempboxb}{\hbox{$\tokend$}}%
  \savebox{\tempboxc}{\hbox{$#6$}}%
\else
  \savebox{\tempboxa}{\hbox{$\hbox{$\tokena$}\times\hbox{$\tokenb$}$}}%
  \savebox{\tempboxb}{\hbox{$\hbox{$\tokena$}\times\hbox{$\tokend$}$}}%
  \savebox{\tempboxc}{\hbox{$\hbox{$\tokena$}\times\hbox{$#6$}$}}%
\fi \ypos=\height \divide\ypos by 2 \xpos=\ypos \advance\xpos by
\width \bfig
\putCtrianglep<-1`1`1;\ypos>(0,0)[`\tokenc`;#5`#6`{#7}]%
\puthmorphism(\ypos,0)[\tokend`\usebox{\tempboxb}`{#8}]{\width}{-1}b%
\puthmorphism(\ypos,\height)[\tokenb`\usebox{\tempboxa}`{#4}]{\width}{-1}a%
\advance\ypos by \width
\putvmorphism(\ypos,\height)[``\usebox{\tempboxc}]{\height}1r%
\efig }}

\def\recurse{\@ifnextchar <{\recursep}{\recursep<\width`\height>}}

\def\puttwohmorphisms(#1,#2)[#3`#4;#5`#6]#7#8#9{{%
\puthmorphism(#1,#2)[#3`#4`]{#7}0a \ypos=#2 \advance\ypos by 20
\puthmorphism(#1,\ypos)[\phantom{#3}`\phantom{#4}`#5]{#7}{#8}a
\advance\ypos by -40
\puthmorphism(#1,\ypos)[\phantom{#3}`\phantom{#4}`#6]{#7}{#9}b }}

\def\puttwovmorphisms(#1,#2)[#3`#4;#5`#6]#7#8#9{{%
\putvmorphism(#1,#2)[#3`#4`]{#7}0a \xpos=#1 \advance\xpos by -20
\putvmorphism(\xpos,#2)[\phantom{#3}`\phantom{#4}`#5]{#7}{#8}l
\advance\xpos by 40
\putvmorphism(\xpos,#2)[\phantom{#3}`\phantom{#4}`#6]{#7}{#9}r }}

\def\puthcoequalizer(#1)[#2`#3`#4;#5`#6`#7]#8#9{{%
\setpos(#1)%
\puttwohmorphisms(\xpos,\ypos)[#2`#3;#5`#6]{#8}11%
\advance\xpos by #8
\puthmorphism(\xpos,\ypos)[\phantom{#3}`#4`#7]{#8}1{#9} }}

\def\putvcoequalizer(#1)[#2`#3`#4;#5`#6`#7]#8#9{{%
\setpos(#1)%
\puttwovmorphisms(\xpos,\ypos)[#2`#3;#5`#6]{#8}11%
\advance\ypos by -#8
\putvmorphism(\xpos,\ypos)[\phantom{#3}`#4`#7]{#8}1{#9} }}

\def\putthreehmorphisms(#1)[#2`#3;#4`#5`#6]#7(#8)#9{{%
\setpos(#1) \settypes(#8)
\if a#9 %
     \vertsize{\tempcounta}{#5}%
     \vertsize{\tempcountb}{#6}%
     \ifnum \tempcounta<\tempcountb \tempcounta=\tempcountb \fi
\else
     \vertsize{\tempcounta}{#4}%
     \vertsize{\tempcountb}{#5}%
     \ifnum \tempcounta<\tempcountb \tempcounta=\tempcountb \fi
\fi \advance \tempcounta by 60
\puthmorphism(\xpos,\ypos)[#2`#3`#5]{#7}{\arrowtypeb}{#9}
\advance\ypos by \tempcounta
\puthmorphism(\xpos,\ypos)[\phantom{#2}`\phantom{#3}`#4]{#7}{\arrowtypea}{#9}
\advance\ypos by -\tempcounta \advance\ypos by -\tempcounta
\puthmorphism(\xpos,\ypos)[\phantom{#2}`\phantom{#3}`#6]{#7}{\arrowtypec}{#9}
}}

\def\setarrowtoks[#1`#2`#3`#4`#5`#6]{%
\def\toka{#1}
\def\tokb{#2}
\def\tokc{#3}
\def\tokd{#4}
\def\toke{#5}
\def\tokf{#6}
}
\def\hex{\@ifnextchar <{\hexp}{\hexp<1000`400>}}
\def\hexp<#1`#2>[#3`#4`#5`#6`#7`#8;#9]{%
\setarrowtoks[#9] \yext=#2 \advance \yext by #2 \xext=#1
\advance\xext by \yext \bfig
\putCtriangle<-1`0`1;#2>(0,0)[`#5`;\tokb``\tokd] \xext=#1
\yext=#2 \advance \yext by #2
\putsquare<1`0`0`1;\xext`\yext>(#2,0)[#3`#4`#7`#8;\toka```\tokf]
\advance \xext by #2
\putDtriangle<0`1`-1;#2>(\xext,0)[`#6`;`\tokc`\toke] \efig }

\makeatother